\documentclass[1p,preprint]{elsarticle}
\usepackage{ecrc_arxiv1}
\volume{00} 
\firstpage{1} 

\makeatletter
\renewcommand\paragraph{\@startsection{paragraph}{4}{\z@}%
           {12\p@ \@plus 6\p@ \@minus 3\p@}%
           {\p@}%
           {\normalfont\normalsize\itshape}}

\makeatother

\usepackage{algorithmic}
\usepackage{amssymb}
\usepackage{amsthm}
\usepackage[title]{appendix}
\usepackage{array}
\usepackage{cancel}
\usepackage{dsfont}
\usepackage{multirow}

\usepackage{epsfig}
\usepackage{calc}
\usepackage{amssymb}
\usepackage{amstext}
\usepackage{amsmath}
\usepackage{hyperref}
\usepackage{lineno}
\usepackage{multicol}
\usepackage{orcidlink}
\usepackage{pslatex}
\usepackage{setspace}
\definecolor{DodgerBlue}{rgb}{0.12, 0.56, 1.0}
\definecolor{Cerise}{rgb}{1, 0.3, 1.0}

\usepackage{xcolor}

\newcommand*{\colorboxed}{}
\def\colorboxed#1#{%
  \colorboxedAux{#1}%
}
\newcommand*{\colorboxedAux}[3]{%
  \begingroup
    \colorlet{cb@saved}{.}%
    \color#1{#2}%
    \boxed{%
      \color{cb@saved}%
      #3%
    }%
  \endgroup
}

\newtheorem{prop}{Proposition}
\newtheorem{thm}{Theorem}
\newtheorem{definition}{Definition}
\begin{document}

\begin{frontmatter}

\title{Accelerating the convergence of Dynamic Iteration method with Restricted Additive Schwarz splitting for the solution of RLC circuits.\tnotemark[1]}

\tnotetext[1]{This work was supported by a grant overseen by the French National Research Agency (ANR) as part of the “Investissements d’Avenir” Program ANE-ITE-002-01.}

\author[supergrid,icj]{Hélèna Shourick \orcidlink{} }
\address[supergrid]{SuperGrid Institute, 23 rue Cyprian,  69200 Villeurbanne, France\\ 
\{helena.shourick,laurent.chedot\}@supergrid-institute.com}

\author[icj]{Damien Tromeur-Dervout\orcidlink{0000-0002-0118-8100}}
\address[icj]{Institut Camille Jordan, Université de Lyon, UMR5208 CNRS-U.Lyon1, Villeurbanne, France\\
	damien.tromeur-dervout@univ-lyon1.fr}
	
\author[supergrid]{Laurent Ch\'edot \orcidlink{} }
\runauth{H. Shourick et al.}

\begin{abstract}
The dynamic iteration method with a restricted additive Schwarz splitting is investigated to co-simulate linear differential algebraic equations system coming from RLC electrical circuit with linear components. We show the pure linear convergence or divergence of the method with respect to the linear operator belonging to the restricted additive Schwarz interface. It allows us to accelerate it toward the true solution with the Aitken's technique for accelerating convergence. This provides a dynamic iteration method less sensitive to the splitting. 
Numerical examples with convergent and divergent splitting show the efficiency of the proposed approach. We also test it on a linear RLC circuit combining different types of circuit modeling (Transient Stability model and Electro-Magnetic Transient model) with overlapping partitions. Finally, some results for a weakly nonlinear differential algebraic equations system are also provided. 
\end{abstract}


\end{frontmatter}


\section{Introduction \label{shourick_contrib_Sec1}}

Since the pioneering work of Lelarasme \& al \cite{lelarasmee1982} that analyze in time domain large-scale problems arising from the modeling of integrated circuits, waveform relaxation methods (WR) \cite{Lumsdaine_White_WR_Semiconductor_NFMAO_1995}  also known as dynamic iteration methods, a term first introduced by Miekkala and Nevanlinna \cite[Eq (2.2)]{Miekkala_Nevanlinna_DI_CVGODE_1987} and generally used in publications \cite{Miekkala_DI_LinearDAE_1989,Arnold_WR_2001, Bartel_Brunk_Guenther_Shoeps_DI_ElectricCircuits_SIAMJSC_2013,  Bartel_Guenther_PDAE_Electrical_SIAMreview_2018, Guenther_Bartel_Jacob_Reis_DI_PPHDAE_CTA_2021} arouses more and more interest with the development of parallel computers
\cite{Reichelt_White_Allen_WRsor_Parallel_SemiCond_SIAMJSC_1995} and more generally in the co-simulation framework \cite{Bartel_Guenther_PDAE_Electrical_SIAMreview_2018,Shoeps_DeGersem_Bartel_Cosimulation_I3ETMagnetics_2012}. 

In such methods applied to Ordinary Differential Equations (ODE) systems or to Differential Algebraic Equations (DAE) systems, the system is decomposed into several subsystems with many internal variables and few external variables. For initial value problems with linear ODEs, the method consists in carrying out some splitting of the linear operator $A=M-N$ \cite{Miekkala_Nevanlinna_DI_CVGODE_1987} such as Jacobi, relaxed Jacobi, Gauss-Seidel or SOR.  Nevertheless this fixed-point process must be  contractant to converge. The analysis of the convergence of the method, using the Laplace transform, says that the convergence  occurs when the spectral radius $\max_{\xi \in \mathbb{R}} \rho( (i\xi I +M)^{-1}N)<1$ \cite[Eq 2.13]{Miekkala_Nevanlinna_DI_CVGODE_1987}. For initial value problems with linear DAE systems  $B \dot{x}+Ax=f$ , like those arising in RLC circuits, Miekkala \cite[Theorem 2]{Miekkala_DI_LinearDAE_1989} extented her previous convergence analysis result with the splitting of $B=M_b-N_b$ and $A=M_A-N_A$, $\max_{\xi \in \mathbb{R}} \rho( (i\xi M_B +M_A)^{-1}(i\xi N_B+N_A))<1$. 
Reichel \& al  combined the waveform SOR with the multistep integration method and showed that the SOR relaxation optimal parameter is dependent on Fourier frequencies \cite[Eq 19]{Reichelt_White_Allen_WRsor_Parallel_SemiCond_SIAMJSC_1995}. Jiang and Wing  determined the expressions of the spectrum and pseudospectrum of the waveform relaxation operators for linear differential-algebraic equations systems which occur especially in circuit simulation \cite[Eqs (3) \& (4)]
{Jiang_Wing_WR_Spectra_LinearDAE_SIAMJNA_2000} and Jiang extended these results to a general class of nonlinear differential-algebraic equations \cite{Jiang_WR_NonlinearDAE_I3ETCircuits_2004} of index one. These extended resuts generalize the expressions of Lumsdaine and Wu\cite{Lumsdaine_Wu_WR_Spectra_SIAMJSC_1997}.
 Several techniques to precondition the fixed point problem were proposed by Arnold \&  Gunther \cite{Arnold_Gunther_preconDI_DAE_Bit_2001}. Hout has established convergence results that are relevant in applications to nonlinear, nonautonomous, stiff initial value problems \cite{Hout_WR_Nonlinear_AppNumMath_1995}.
 
Some convergence acceleration techniques for the WR have been proposed. Some waveform successive overrelaxation (SOR) techniques have been proposed by Janssen and Vandewalle  \cite{Janssen_Vandewalle_WR_SOR_SIAMJNA_1997} to accelerate the standard waveform method. Leimkuhler proposed  to accelerate the WR by solving the defect equations with a larger timestep, or by using a recursive procedure based on a succession of increasing timesteps \cite{Leimkuhler_WR_TimeStepAcceleration_SIAMJNA_1998}. Lumdaisne \& Wu proposed to accelerate the WR by  Krylov subspace techniques (WGMRES) \cite{Lumsdaine_Wu_WR_KrylovAcceleration_SIAMJNA_2003} to solve time-dependent problems. Botchev \& al \cite{Botchev_Oseledets_Tyrtyshnikov_WR_Krylov_ComputMathAppl_2014} compared WR-Krylov with Krylov's methods combined with the shift and invert (SAI) technique to obtain parallelism in time. Ladics \cite{Ladics_WR_SemilinearPDE_ErrorEstimate_ComputApplMath_2015} combined the WR with convergent numerical methods to solve semi-linear PDEs, he showed the effect of  applying time windows. Recent developments in the dynamic iteration method for the co-simulation of electrical circuits have been carried out by Bartel \& al  \cite{Bartel_Brunk_Schoeps_DI_cvgrate_JCAppMath_2014,Bartel_Brunk_Guenther_Shoeps_DI_ElectricCircuits_SIAMJSC_2013} and by Ali \&  al \cite{ Ali_Bartel_Brunk_Schoeps_DI_cvg_ProcMI_2012}. Gausling \& al \cite{Gausling_Bartel_DI_COSIM_CG_Stoch_ProcMI_2018}  analyzed the contraction and the rate of convergence of the co-simulation process for a test circuit subjected to uncertainties on the parameters of its components. Morever, the rate of convergence or divergence of the dynamic iteration depends on the interface coupling \cite{Gausling_Bartel_DI_Cosim_cvgInterface_2016}. Pade and Tischendorf \cite{Pade_Tischendorf_WR_CVGCriterion_NUMAlgo_2019} presented topological criteria for the coupling of networks which are easy to check and which are sufficient to ensure the convergence of the WR which is related to the DAE index.

In this paper, we focus on dynamic iteration for  linear DAE as in\cite{Jiang_Wing_WR_Spectra_LinearDAE_SIAMJNA_2000} with a perspective of accelerating domain decomposition. Garbey and Tromeur Dervout developed the  Aitken-Schwarz domain decomposition \cite{tromeur_contrib_DTD} using the pure linear convergence of Schwarz type method to accelerate their convergence with the Aitken's technique for accelerating convergence for PDE  problems. Tromeur-Dervout \cite{tromeur_contrib_ISI:000268885300011,tromeur_contrib_DTD2} developed a completely algebraic formulation of this Aitken's technique for accelerating convergence as a fonction of the trace of the Schwarz's iterations on the domain decomposition interface. Shourick \& al developed heterogenous Schwarz domain decomposition accelerated by the Aitken's technique for accelerating convergence for the co-simulation of electromagnetic transient and transient stability \cite{Shourick_DD26}. We want to demonstrate with such a technique that we can  acceleterate the dynamic iteration to obtain the solution whatever its convergence or its divergence and therefore without topological criteria on the coupling as in \cite{Pade_Tischendorf_WR_CVGCriterion_NUMAlgo_2019}   or \cite{Gausling_Bartel_DI_Cosim_cvgInterface_2016}.

The outline of this paper is as follows: section \ref{linearDAE} focuses on dynamic iteration for linear DAEs for which the splitting of linear operators follows a Restrictive Additive Schwarz (RAS) domain decomposition. Section \ref{Dynamic_iteration_DDM} details the error operator of the dynamic iteration associated with the RAS domain decomposition and the acceleration of the convergence with Aitken's technique using the pure linear convergence, i.e. the error operator does not depend on the iteration. We present  numerical results  in section \ref{result} with examples of linear DAEs of index one and two, as well as with heterogeneous modeling using ElectroMagnetic Transient modeling and Transient Stability modeling with overlap. We conclude in section \ref{conclusion}.

\section{Dynamic iteration for linear DAE \label{linearDAE}}
Let us consider the linear DAE
\begin{eqnarray}
\left\{\begin{array}{rcl} \dot{x}(t) +A x(t)+B y(t) &=& b_1(t),  x(0)=x_0,   \\
 C x(t) + D y(t)  &=& b_2(t),\, t\in [0,T]. \end{array} \right. \label{EqLinearDAE}
\end{eqnarray}

where $x(t)\in \mathbb{R}^{n_1}$ and $y(t) \in  \mathbb{R}^{n_2}$ for all $t \in [0, T]$, $D$  is a $n_2 \times n_2$ nonsingular matrix, $A$ is an $n_1 \times n_1$ matrix, $B$ is an $n_1 \times n_2$ matrix,
$C$ is an $n_2 \times n_1$ matrix, $b_1(t) \in \mathbb{R}^{n_1}$ and $b_2(t) \in \mathbb{R}^{n_2}$ are known input functions, and
$x_0 \in \mathbb{R}^{n_1}$ is a consistent initial value. Let $n = n_1 + n_2$.

Let us write the dynamic iteration in the context of a RAS domain decomposition. First, we define the matrix $\mathbb{A}=\left(\begin{array}{cc} A & B \\ C & D \end{array} \right)$ corresponding to the linear operator of the DAE and we define $z(t)=[x(t),y(t)]^T$, $b(t)=[b_1(t), b_2(t)]^T$ and $\mathbb{I}_{d}=\left(\begin{array}{cc} I_{n_1} & 0_{n_1\times n_2} \\ 0_{n_2 \times n_1} & 0_{n_2\times n_2} \end{array} \right)$. Then we can rewrite Eq \eqref{EqLinearDAE} as:
\begin{eqnarray}
\begin{array}{rcl} \mathbb{I}_{d}\dot{z}(t) + \mathbb{A} z(t) &=& b(t),\;  x(0)=x_0, \; t\in [0,T]. \end{array}  \label{EqLinearDAE2}
\end{eqnarray}
By adapting the notations of \cite{tromeur_contrib_RAS}, we consider the matrix $\mathbb{A}\in \mathbb{R}^{n\times n}$ having a non-zero pattern and the associated graph $G = (W, F)$, where the set of vertices $W = \left\{1,\ldots, n\right\}$ represents the $n$ unknowns and the set of edges $F = \left\{(i, j) | \mathbb{A}_{i,j} \neq 0\right\}$ represents the pairs of vertices that are coupled by a non-zero element in $\mathbb{A}$. Then, we assume that a graph partitioning was applied and that resulted in $N$ non-overlapping subsets $W_i^0$ whose union is $W$. 
Let $W_i^p$ be the $p$-overlap partition of $W$, obtained by including all the vertices immediately neighboring the vertices of $W_i^{p-1}$. Let $W_{i,e}^p=W_i^{p+1} \backslash W_i^{p}$. Then let $R_i^p \in \mathbb{R}^{n_i \times n}$  (  $R_{i,e}^p \in \mathbb{R}^{n_{i,e} \times n}$ and $\tilde{R}_i^0 \in \mathbb{R}^{n_i \times n}$ respectively) be the operator which restricts $x\in \mathbb{R}^n$ to the components of $x$ belonging to $W_i^p$   ($W_{i,e}^p$  and $W_i^0$ respectively, and the operator $\tilde{R}_i^0 \in \mathbb{R}^{n_i \times n}$ puts $0$ to the unknowns belonging to $W_i^p\backslash W_i^0$). 
 
We define the operators $\mathbb{A}_i = R_i^p \mathbb{A} R_i^{pT}$ and $\mathbb{E}_i=R_i^p \mathbb{A} R_{i,e}^{pT}$, the vectors $z_i = R_i^p z$, $b_i = R_i^p b$, and $z_{i,e} =R_{i,e}^p z$. We also introduce $R_{i}^{p,d} \in \mathbb{R}^{n_{i_1} \times n_1}$ and $R_{i}^{p,a}\in \mathbb{R}^{n_{i_2} \times n_2} $ (respectively $\tilde{R}_{i}^{p,d} \in \mathbb{R}^{n_{i_1} \times n_1}$, $\tilde{R}_{i}^{0,a}\in \mathbb{R}^{n_{i_2} \times n_2} $  $R_{i,e}^{p,d}\in \mathbb{R}^{n_{{ie}_1} \times n_1}$ and $R_{i,e}^{p,a}\in \mathbb{R}^{n_{{ie}_2} \times n_2}$) which separate the differential and the algebraic variables belonging to $W_i^p$ (respectively $W_i^0$ and $W_{i,e}^p$) i.e.
$R_i^{p}= \left( \begin{array}{cc} R_i^{p,d} & 0_{n_{i_1}\times n_2}  \\ 0_{n_{i_2}\times n_1} &  R_i^{p,a} \end{array} \right)$ (respectively $\tilde{R}_i^{0}= \left( \begin{array}{cc} \tilde{R}_i^{0,d} & 0_{n_{i_1}\times n_2}  \\ 0_{n_{i_2}\times n_1} &  \tilde{R}_i^{0,a} \end{array} \right)$  and $R_{i,e}^{p}= \left( \begin{array}{cc} R_{i,e}^{p,d} & 0_{n_{{ie}_1}\times n_2}  \\ 0_{n_{{ie}_2}\times n_1} &  R_{i,e}^{p,a} \end{array} \right)$) and\\
 $\mathbb{I}_{i,d} = \left( \begin{array}{cc} I_{n_{i_1}} & 0_{n_{{i}_1}\times {n_{i_2}}}  \\ 0_{n_{{i}_2}\times n_{i_1}} &  0_{n_{{i}_2}\times n_{{i}_2}} \end{array} \right)$.

\begin{definition}
The dynamic iteration ${k+1}$ with the RAS splitting is written locally for the $W_i^p$ partition as:
\begin{eqnarray}
\left\{\begin{array}{rcl} \mathbb{I}_{i,d} \dot{z_i}^{(k+1)}(t) + \mathbb{A}_i z_i^{(k+1)}(t) &=& b_i(t)-\mathbb{E}_i z_{i,e}^{(k)}(t),\; t\in [0,T],\\
  x_i(0)=  R_{i}^{p,d} x_0. \end{array} \right. \label{EqLinearDAERAS}
\end{eqnarray}
and by separating the differential and algebraic variables belonging to the partition $W_i^p$ :
\begin{eqnarray}
\left\{\begin{array}{rcl} \dot{x}_i^{(k+1)}(t) +A_i x_i^{(k+1)}(t)+B_i y_i^{(k+1)}(t) &=& b_{i,d} (t)- E_{i,d}^{d} x_{ie}^{(k)}(t)  - E_{i,d}^{a} y_{ie}^{(k)}(t) ,   \\
 C_i x_i^{(k+1)}(t) + D_i y_i^{(k+1)}(t)  &=& b_{i,a}(t)- E_{i,a}^{d} x_{ie}^{(k)}(t)  - E_{i,a}^{a} y_{ie}^{(k)}(t),\\
  x_i(0)=R_i^{p,d} x_0,&& \, \, t\in [0,T]. \end{array} \right. \label{EqLinearDAERASsplit}
\end{eqnarray}
Where 
\begin{eqnarray}
\left( \begin{array}{cc} A_i & B_i \\ C_i & D_ i \end{array} \right) &=& \left( \begin{array}{cc} R_i^{p,d} & 0_{n_{i_1} \times n_{2}} \\  0_{n_{i_2} \times n_{1}} & R_i^{p,a} \end{array} \right) \mathbb{A} \left( \begin{array}{cc}  (R_i^{p,d})^T & 0_{ n_{1} \times n_{i_2}} \\   0_{n_{2}\times n_{i_1} } & (R_i^{p,a})^T  \end{array} \right),\\
\left( \begin{array}{cc} E_{i,d}^d &E_{i,d}^a \\ E_{i,a}^d &E_{i,a}^a \end{array} \right) &=& \left( \begin{array}{cc} R_i^{p,d} & 0_{n_{i_1} \times n_{2}} \\  0_{n_{i_2} \times n_{1}} & R_i^{p,a} \end{array} \right) \mathbb{A} \left( \begin{array}{cc}  (R_{i,e}^{p,d})^T & 0_{ n_{1} \times n_{{ie}_2}} \\   0_{n_{2}\times n_{{ie}_1} } & (R_{i,e}^{p,a})^T  \end{array} \right).
\end{eqnarray}
\end{definition}

\begin{prop}
The dynamic iteration with the RAS splitting is written globally as:
\begin{eqnarray}
\left\{\begin{array}{rcl} 
\dot{x}^{(k+1)}(t) + A_1^{d}  x^{(k+1)}(t)+ B_1^{d} y^{(k+1)}(t)&=&  b^{d} (t)+ A_2^{d} x^{(k)}(t)  + B_2^{d}  y^{(k)}(t) ,   \\
C_1^{a} x^{(k+1)}(t) + D_1^{a} y^{(k+1)}(t)  &=& b^{a}(t)+ C_2^a x^{(k)}(t)  + D_2^a y^{(k)}(t),\\
  x^{(k+1)}(0)=x_0,&& \, \, t\in [0,T]. \end{array} \right.
\end{eqnarray}
with
\begin{eqnarray}
A_1^{d} &=&\sum_{i=0}^{N-1} \tilde{R}^{0,d}_i A_i R^{p,d}_i, \quad
A_2^{d} = -\sum_{i=0}^{N-1} \tilde{R}^{0,d}_i E_{i,d}^d R^{p,d}_{ie}, \quad  b^{d}(t)=\sum_{i=0}^{N-1} \tilde{R}_i^{0,d} R^{p,d}_i b (t), \nonumber\\
B_1^{d}&=&\sum_{i=0}^{N-1} \tilde{R}^{0,d}_i B_i R^{p,a}_i, \quad
B_2^{d} =-\sum_{i=0}^{N-1} \tilde{R}^{0,d}_i E_{i,d}^a R^{p,a}_{ie}, \nonumber \\
C_1^{a} &=&\sum_{i=0}^{N-1} \tilde{R}^{0,a}_i C_i R^{p,d}_i, \quad
C_2^{a} =-\sum_{i=0}^{N-1} \tilde{R}^{0,a}_i E_{i,a}^d R^{p,d}_{ie},   \nonumber\\
D_1^{a} &=&\sum_{i=0}^{N-1} \tilde{R}^{0,a}_i D_i R^{p,a}_i, \quad
D_2^{a} =-\sum_{i=0}^{N-1} \tilde{R}^{0,a}_i E_{i,a}^a R^{p,a}_{ie},\quad b^{a}(t) = \sum_{i=0}^{N-1} \tilde{R}_i^{0,d} R^{p,a}_i b (t). \nonumber
\end{eqnarray}
\begin{proof}
The sum of the contribution of each partition with $\sum_{i=0}^{N-1} \tilde{R}^{0,d}_i $ and the definitions of $x_i$, $y_i$,$x_{i,e}$, $y_{i,e}$ give the result.
\end{proof}
\end{prop}

\section{Dynamic Iteration error operator and acceleration\label{Dynamic_iteration_DDM}}
We are in the formalism of the dynamic iteration with a general splitting. By adapting the results \cite[Eq. (3) and  (4)]{Jiang_Wing_WR_Spectra_LinearDAE_SIAMJNA_2000} to our notations, we have:
\begin{thm}[Jiang \& Wing \cite{Jiang_Wing_WR_Spectra_LinearDAE_SIAMJNA_2000}]
The dynamic iteration with RAS splitting applied to a linear DAE system has an error operator $\mathcal{R}$ which does not depend on the iteration number such as:
\begin{eqnarray}
z^{(k)} &=& \mathcal{R} z^{(k-1)} + \phi \label{DI_linear_z}\\
\mathcal{R} &=& \left( \begin{array}{cc}\mathcal{R}_1 & \mathcal{R}_2 \label{DI_linear_R},\\
(D_1^a)^{-1}( C_2^d- C_1^d \mathcal{R}_1)& (D_1^a)^{-1} (D_2^a-C_1^d \mathcal{R}_2)\end{array} \right).
\end{eqnarray}
With 
\begin{eqnarray}
(\mathcal{R}_1u)(t) &=& \int_0^t e^{-S_1(t-s)} (A_2^d-B_1^d(D_1^a)^{-1}C_2^a) u(s) ds, \forall  u \in L^2([0,T], \mathbb{C}^{n_1}),\\
(\mathcal{R}_2v)(t) &=& \int_0^t e^{-S_1(t-s)} (B_2^d-B_1^d(D_1^a)^{-1}D_2^a) v(s) ds ,\forall  v \in L^2([0,T], \mathbb{C}^{n_2}),\\
S_1&=& A_1^d-B_1^d(D_1^a)^{-1} C_1^a, \, \phi(t)=[\phi_1(t),\phi_2(t)],\\
\phi_1(t)&=& e^{-S_1 t} x_0+ \int_0^t  e^{-S_1(t-s)}[b_1^d(t)-B_1^d (D_1^a)^{-1} b_2^a(t)]ds, \\
\phi_2(t)&=&  -(D_1^a)^{-1} (C_1^d \phi_1(t)+ b_2^a(t)). 
\end{eqnarray}
\end{thm}

The interest of Eqs \eqref{DI_linear_z} \& \eqref{DI_linear_R} is to show the pure linear convergence of the DI and the possibility of accelerating the convergence 
to the true solution $z^{(\infty)}$ with the Aitken's technique for accelerating convergence, if 1 is not an eigen value of $\mathcal{R}$, as follows:
\begin{eqnarray}
z^{(\infty)} &=&  (I-\mathcal{R})^{-1} (z^{(1)} + \mathcal{R} z^{(0)})
\end{eqnarray}
We present now the discrete counterpart of the DI with  RAS splitting and its Aitken's technique for accelerating convergence  involving the interface solution of the RAS. We use a backward Euler for time discretization, other backward differences formula (BDF) schemes would give similar results with more complicated formula. 

\begin{eqnarray}
\left\{\begin{array}{rcl} 
 \tilde{A}_1^{d}  x^{n+1,(k+1)}+  \tilde{B}_1^{d} y^{n+1,(k+1)}&=&   \tilde{b}^{n+1,d} +  \tilde{A}_2^{d} x^{n+1,(k)}  +  \tilde{B}_2^{d}  y^{n+1,(k)} ,   \\
C_1^{a} x^{n+1,(k+1)} + D_1^{a} y^{n+1,(k+1)}  &=& b^{n+1,a}+ C_2^a x^{n+1,(k)} + D_2^a y^{n+1,(k)},\\
  x^{0,(k+1)}=x_0.&&  \end{array} \right.
\end{eqnarray}\label{DAEDiscret}
with $\tilde{A}_1^{d}= I+\Delta t \; {A}_1^{d}$, $ \tilde{B}_1^{d}=\Delta t \; {B}_1^{d}$, $ \tilde{A}_2^{d}= \Delta t \; {A}_2^{d}$, $ \tilde{B}_2^{d}=\Delta t \; {B}_2^{d}$, $ \tilde{b}^{n+1,d} = x^{n,*}+ \Delta t \;{b}^{n+1,d}$, with $x^{n,*}=x^{n,k+1}$ or $x^{n,*}=x^{n,\infty} $ depending on the implementation strategy used in section \ref{strategy}. 

Locally, it is written with $x_i^{0,({k+1})}=R_i^{p,d} x_0$:
\begin{eqnarray}
\underbrace{\left(\begin{array}{c}  x_i^{n+1,(k+1)}\\ y_i^{n+1,(k+1)}\end{array} \right)}_{z_i^{n+1,(k+1)}} &=& \underbrace{\left( \begin{array}{cc} \tilde{A}_i & \tilde{B}_i \\  C_i & D_i \end{array} \right)^{-1}}_{\tilde{\mathbb{A}}_i^{-1}}   \underbrace{ \left( \left(\begin{array}{c} \tilde{b}_{i,d}^{n+1} \\  b_{i,a}^{n+1}
\end{array} \right) \right.}_{\tilde{b}^{n+1}_i} - \underbrace{\left( \begin{array}{cc} \tilde{E}_{i,d}^{d} & \tilde{E}_{i,d}^{a} \\ E_{i,a}^{d} &  E_{i,a}^{a} \end{array} \right)}_{\tilde{\mathbb{E}}_i} \underbrace{\left. \left( \begin{array}{c} x_{i,e}^{n+1,(k)}  \\ y_{i,e}^{n+1,(k)} \end{array} \right) \right)}_ {z_{i,e}^{n+1,(k)}}   \label{EqLinearDAERASsplitDiscret2}
\end{eqnarray}
 
 By defining $M_{RAS}^{-1}\stackrel{def}{=} \sum_{i=0}^{N-1} \tilde{R}_i^{0T} \tilde{\mathbb{A}}_i^{-1} R_i^p$ and adding the contribution of each partition $W_i^p$, the RAS can be seen as a Richardson's process:
 
\begin{eqnarray}
\sum_{i=0}^{N-1} \tilde{R}_i^{0T} R_i^p z^{n+1,(k+1)} &=&\sum_{i=0}^{N-1} \tilde{R}_i^{0T} \tilde{\mathbb{A}}_i^{-1} R_i^p \tilde{b}^{n+1} -\sum_{i=0}^{N-1} \tilde{R}_i^{0T} \tilde{\mathbb{A}}_i^{-1} R_i^p \tilde{\mathbb{A}} R_{i,e}^{pT} R_{i,e}^{p} z^{n+1,(k)} \label{TromeurDervout_contrib_eqRAS2},\\
z^{n+1,(k+1)} &= & M_{RAS}^{-1} \tilde{b}^{n+1} - M_{RAS}^{-1} \tilde{\mathbb{A}}z^{n+1,(k)} + z^{n+1,(k)}, \nonumber\\
&=& z^{n+1,(k)}+ M_{RAS}^{-1}(\tilde{b}^{n+1}-\tilde{\mathbb{A}}z^{n+1,(k)}). \label{TromeurDervout_contrib_equation:aitken:richardson}
\end{eqnarray}
The Richardson's process \eqref{TromeurDervout_contrib_equation:aitken:richardson} is deduced from \eqref{TromeurDervout_contrib_eqRAS2} (see \cite[Theorem 3.7]{tromeur_contrib_Gander08schwarzmethods}) by using the property  $R_i^p\tilde{\mathbb{A}}=R_i^p\tilde{\mathbb{A}}(R_i^{pT}R_i^p  +R_{i,e}^{pT} R_{i,e}^p)$.

The restriction of \eqref{TromeurDervout_contrib_equation:aitken:richardson} to the interface $\Gamma=\left\{W_{0,e}^p,\ldots, W_{N-1,e}^p \right\}$ of size\\ $n_\Gamma=\sum_{i=0}^{N-1} n_{i,e}$, by setting $R_\Gamma = (R_{0,e}^p,\ldots, R_{N-1,e}^p)^T \in \mathbb{R}^{n_\Gamma \times n}$ and using the property\\ $R_{i,e}^{pT} R_{i,e}^p R_\Gamma^T R_\Gamma=R_{i,e}^{pT} R_{i,e}^p$, can be written as: 
\begin{equation}
\begin{array}{lcl}
\underbrace{R_\Gamma  z^{n+1,(k+1)}}_{z_\Gamma^{n+1,(k+1)}} 
        &=&\underbrace{ R_\Gamma \left(I - M^{-1}_{RAS} \tilde{\mathbb{A}} \right)R_\Gamma^T}_{P^{~}} \underbrace{R_\Gamma z^{n+1,(k)}}_{z_\Gamma^{n+1,(k)}}+\underbrace{R_\Gamma M^{-1}_{RAS} \tilde{b}^{n+1}}_{c^{n+1}}.
       
\label{TromeurDervout_contrib_eq:rich}
\end{array}
\end{equation}
The pure linear convergence of the RAS at the interface given by: $z_\Gamma^{n+1,(k)}-z_\Gamma^{n+1,(\infty)}=P (z^{n+1,{(k-1)}}-z_\Gamma^{n+1,(\infty)})$ (the error operator $P$ does not depend on the iteration $k$) allows to apply the Aitken's technique for accelerating convergence  to obtain the true solution  $z_\Gamma^{n+1,(\infty)}$ on the interface $\Gamma$: $z_\Gamma^{n+1,(\infty)} = (I-P)^{-1} (z_\Gamma^{n+1,(k)}-Pz_\Gamma^{n+1,(k-1)})$, and thus after another local resolution, the true solution $z^{n+1,(\infty)}$. Let us note that one can accelerate the convergence toward the solution for an iterative convergent or divergent method. The only need is that 1 is not one of the  eigen values of $P$. Considering $e^k=z_\Gamma^{n+1,(k)}-z_\Gamma^{n+1,(k-1)}, k=1,\ldots $, the  operator $P \in \mathbb{R}^{n_\Gamma \times n_\Gamma}$ can be computed algebraically after $n_\Gamma+1$ iterations as $P=[e^{k_\Gamma+1},\ldots, e^2][e^{k_\Gamma},\ldots, e^1]^{-1}$.
Let us notice the sparse structure of the operator $P$ in the two partitions case as we will have in the numerical examples. Defining $e_i$ the restriction of the error $e$ to the partition $W_{i,e}^{p}$ leads to write the DI with RAS splitting as:

\begin{eqnarray}
\left(
\begin{array}{c}
e_{0}\\
e_{1}\\
\end{array}  \right) ^{k+1}=
\left( \begin{array}{cc}
0& P_{0} \\
P_{1} & 0
\end{array}  \right) 
\left(
\begin{array}{c}
 e_{0}\\
e_{1}
\end{array}  \right) ^{k}
\end{eqnarray}
 Using the properties of $\sum_{i=0}^{N-1} \tilde{R}_i^0$ and $R_{i,e}^{p}$, $P_i$ is given by :
 \begin{eqnarray}
 P_0 &=&  R_{0,e}^{p} \tilde{\mathbb{A}}_0^{-1} R_0^p \tilde{\mathbb{A}} R_{1,e}^{pT} \label{contribCalcul0}\\
 P_1 &=&  R_{1,e}^{p} \tilde{\mathbb{A}}_1^{-1} R_1^p \tilde{\mathbb{A}} R_{0,e}^{pT}
 \label{contribCalcul1}
 \end{eqnarray}
 
 \section{Strategies for Dynamic Iteration}\label{strategy}
 
 \subsection{Sequential time steps strategy}
 In the sequential time step strategy, we apply the Aitken's technique for accelerating convergence, after $n_\Gamma+1$ DI iterations for the first regular time step, in order to numerically build the $P$ operator. Then, if we use the same time step size for the following time steps, we can perform the Aitken's convergence acceleration technique after one DI iteration.
 
 \subsection{Pipelined time steps strategy}
In the pipelined time step strategy, we perform several time steps per DI iteration. Then after $n_\Gamma+1$ we can compute the $P$ operator associated with a time step, and we can build the $\mathbb{P}$ operator error associated with these several time steps given in \eqref{Aitken_Several}.
 \begin{eqnarray}
 \underbrace{\left( 
 \begin{array}{c}  
 z^{1,(k+1)}_\Gamma \\
 z^{2,(k+1)}_\Gamma \\
 \vdots\\
 z^{m-1,(k+1)}_\Gamma \\
 z^{m,(k+1)}_\Gamma 
 \end{array}\right)}_{Z_\Gamma^{(k+1)}} &=& \underbrace{\left( 
\begin{array}{ccccc} 
 P &&&&\\  
  \mathbb{I} &P &&&\\ 
  &\ddots & \ddots &&\\
  &&\mathbb{I} &P &\\
  &&&\mathbb{I} &P 
 \end{array}\right)}_{\mathbb{P}} \underbrace{\left( 
 \begin{array}{c} 
 z^{1,(k)}_\Gamma \\
 z^{2,(k)}_\Gamma \\
 \vdots\\
 z^{m-1,(k)}_\Gamma \\
 z^{m,(k)}_\Gamma 
 \end{array}\right)}_{Z_\Gamma^{(k)}} + \underbrace{\left( 
 \begin{array}{c} 
 z_0+c_1 \\
 c_2 \\
 \vdots\\
 c_{m-1} \\
c_m \\
 \end{array}\right)}_{C} \label{Aitken_Several}
 \end{eqnarray}
 
 This pipelined time steps strategy has also a pure linear convergence/divergence and can also be accelerated by the Aitken's technique for accelerating convergence .

\section{Numerical results for DI with the RAS splitting\label{result} }

Firstly, we test our method on the  RLC examples of \cite{Pade_Tischendorf_WR_CVGCriterion_NUMAlgo_2019}  that they use to illustrate their convergence criteria for the WR Gauss-Seidel  (the first example converges and the second, an index 2 DAE, diverges). We will see that the convergence or divergence of the method depends on the time step chosen for the DI with the RAS splitting. Finally, we apply the  Aitken's technique for accelerating convergence   on the convergent and divergent cases. It shows the  possibility of the method to obtain the true solution even in the divergent cases.

\subsection{First example of \cite{Pade_Tischendorf_WR_CVGCriterion_NUMAlgo_2019}}
The first RLC circuit example, satisfies the  criteria of  \cite{Pade_Tischendorf_WR_CVGCriterion_NUMAlgo_2019} to ensure the convergence for the WR Gauss-Seidel method. The circuit splitting is as follows:

\begin{minipage}{10.5cm}
\begin{minipage}{5.2cm}
 \begin{tikzpicture}[scale=0.55]
 \draw[Cerise]  (-6,2) -- (-6,1.15); 
  \draw[DodgerBlue]  (-6,1) -- (-6,-2); 
 
 \draw[Cerise]  (-6,2) -- (-4.5,2);
  \draw[Cerise]  (-3.5,2) -- (-2.5,2);
   \draw[Cerise]  (-2.5,2) -- (-2.5,1.15);

  \draw[DodgerBlue]  (1,1) -- (1,-2); 
 
  \draw[DodgerBlue]    (1,1)--(-0.5,1); 
    \draw[DodgerBlue]  (-1.5,1) -- (-3.5,1); 
   \draw[DodgerBlue]  (-4.5,1) -- (-6,1); 
   
    \draw[DodgerBlue]  (-2.5,1) -- (-2.5,-0.3); 
     \draw[DodgerBlue]  (-2.5,-0.7) -- (-2.5,-2); 
     
       \draw[DodgerBlue]  (1,-2)--(-0.5,-2); 
    \draw[DodgerBlue]  (-1.5,-2)-- (-3.5,-2); 
   \draw[DodgerBlue]  (-4.5,-2) -- (-6,-2); 
     
  \fill[DodgerBlue] (-6,1) circle(0.07);
   \fill[DodgerBlue] (1,1) circle(0.07);
  \fill[DodgerBlue] (-2.5,-2) circle(0.07);
   \fill[DodgerBlue] (-2.5,1) circle(0.07);
   
  \draw[DodgerBlue] (-6,1) node[left]{\scriptsize $n_1$};
   \draw[DodgerBlue] (1,1) node[right]{\scriptsize $n_r$};
   \draw[DodgerBlue]  (-2.5,-2) node[below]{\scriptsize $n_3$};
   \draw[DodgerBlue] (-2.7,1) node[above]{\scriptsize $n_2$};
   
   \begin{scope}[shift={(-5,2)},rotate=90]
{
 \draw[white, fill=white] (-0.1,0) rectangle (0.4,-1.5);
 \foreach \r in {0,...,2}
 {
  \draw[Cerise,thick,scale=1/3,shift={(0,-\r)}]
	(0,0) .. controls ++(2,0) and ++(1,0) ..
	++(0,-1.5) .. controls ++(-1,0) and ++(-0.5,0) ..
	++(0,0.5);
 }
 \draw[Cerise,thick,scale=1/3,shift={(0,-3)}] (0,0) .. controls ++(2,0) and ++(1,0) .. ++(0,-1.5);
}

\end{scope}
\draw[Cerise]  (-4,2.4) node[above]{\footnotesize  $L_1$};

        \begin{scope}[shift={(-1.5,1)},rotate=90]
{
 \draw[white, fill=white] (-0.1,0) rectangle (0.4,-1.5);
 \foreach \r in {0,...,2}
 {
  \draw[DodgerBlue,thick,scale=1/3,shift={(0,-\r)}]
	(0,0) .. controls ++(2,0) and ++(1,0) ..
	++(0,-1.5) .. controls ++(-1,0) and ++(-0.5,0) ..
	++(0,0.5);
 }
 \draw[DodgerBlue,thick,scale=1/3,shift={(0,-3)}] (0,0) .. controls ++(2,0) and ++(1,0) .. ++(0,-1.5);
}
\end{scope}
 \draw[DodgerBlue](-0.6,0.21) node[above]{\footnotesize $L_2$};

   \draw[DodgerBlue,ultra thick] (-2.1,-0.3) -- (-2.8,-0.3);
   \draw[DodgerBlue,ultra thick] (-2.1,-0.7) -- (-2.8,-0.7);

  \draw[DodgerBlue] (-1.8,-0.5) node{\footnotesize  $C$};

   \draw[DodgerBlue,thick] (-4,-2) circle(0.5);
  \draw[DodgerBlue,thick] (-3.5,-2)-- (-4.5,-2);

  \draw[DodgerBlue] (-4,-2.8) node{\footnotesize $E_v $};

   \draw[DodgerBlue,thick] (-1,-2) circle(0.5);
  \draw[DodgerBlue,thick] (-1,-2.5)-- (-1,-1.5);

  \draw[DodgerBlue] (-1,-2.8) node{\footnotesize $E_i $};

    \draw[DodgerBlue](-4.5,1.15) -- (-3.5,1.15); 
   \draw[DodgerBlue] (-4.5,0.85) -- (-3.5,0.85); 
    \draw[DodgerBlue](-4.5,1.15) -- (-4.5,0.85); 
   \draw [DodgerBlue](-3.5,1.15) -- (-3.5,0.85);

   \draw[DodgerBlue] (-4,0.21) node[above]{\footnotesize $G$};
 \end{tikzpicture}
\vspace*{-0.5cm}
\begin{eqnarray}
\small
\color{Cerise} e_1-e_2-L_1 \dfrac{d{i}_{1}}{dt} & \color{Cerise}=& \color{Cerise}0, \nonumber \\
\color{DodgerBlue} e_r&\color{DodgerBlue} =&\color{DodgerBlue} 0, \nonumber\\
\color{DodgerBlue}e_1-e_3-E_{v}(t)-Z_s i_{5} &\color{DodgerBlue} =& \color{DodgerBlue} 0,\nonumber 
\end{eqnarray}
\end{minipage}
\hfill
\begin{minipage}{5.4cm}
\begin{eqnarray}
\small
\color{DodgerBlue}G (e_1-e_2) -i_{2} &\color{DodgerBlue} =& \color{DodgerBlue} 0,\nonumber \\
\color{DodgerBlue}e_2-e_r-L_2 \dfrac{d{i}_{3}}{dt} &\color{DodgerBlue} =&\color{DodgerBlue}  0, \nonumber\\
\color{DodgerBlue}C (\dfrac{d{e}_3}{dt}-\dfrac{d{e}_2}{dt})-i_{4} &\color{DodgerBlue} =&\color{DodgerBlue} 0,\nonumber\\
\color{DodgerBlue}E_{i}(t)-i_{6} &\color{DodgerBlue} =& \color{DodgerBlue} 0,\nonumber \\
\color{DodgerBlue}-i_{1}-i_{2}+i_{3}-i_{4}&\color{DodgerBlue} =&\color{DodgerBlue}  0, \nonumber \\
\color{DodgerBlue}i_{1}+i_{2}+i_{5}&\color{DodgerBlue} =&\color{DodgerBlue}  0, \nonumber \\
\color{DodgerBlue}-i_{5}+i_{4}+i_{6}&\color{DodgerBlue} =&\color{DodgerBlue}  0.  \nonumber
\end{eqnarray}
\end{minipage}
\end{minipage}
\newline
\medskip

To adapt this circuit to the formalism of Eq.\eqref{EqLinearDAE}, it is necessary to define new variables which are the combination of variables whose derivatives are combined in the same equation and to add the corresponding algebraic combination to the equations of the system. Here, ${ \color{DodgerBlue}C (\dfrac{dv_1}{dt})-i_{4} =0}$ and ${\color{DodgerBlue}v_1=e_3-e_2}$ where $v_1$ is part of $ x$ and $e_3,e_2$ are part of $y$.

The interface values are  $e_1$ and $e_2$ for the first partition and $i_1$ for the second partition.
The error operator ${P}$ and it's eigen values are calculated following Eq.  \eqref{contribCalcul0} and \eqref{contribCalcul1}:
\medskip
\newline
\begin{minipage}{14cm}
\begin{minipage}{8cm}
\begin{eqnarray}
 {P}&=&\left(\begin{array}{ccc}
0 &\colorboxed[rgb]{0.12, 0.56, 1.0}{\begin{matrix}  -\frac{\Delta t}{L_1} &\frac{\Delta t}{L_1} \end{matrix}}\\
\colorboxed[rgb]{1, 0.3, 1.0}{-\frac{L_2}{\Delta t}+\frac{\frac{L_2C}{\Delta t^2}+\frac{L_2G}{\Delta t}+1}{\frac{C}{\Delta t}+G}}&0&0\\
0&0&0
\end{array}  \right),  \nonumber
\end{eqnarray}
\end{minipage}
\hfill
\begin{minipage}{4cm}
\begin{tabular}{c|c}
$\lambda_1$ & 0\\
\hline
$\lambda_2$ & $-i \sqrt{\frac{\Delta t}{L_1(\frac{C}{\Delta t}+G)} }$\\
\hline
$\lambda_3$ & $i \sqrt{\frac{\Delta t}{L_1(\frac{C}{\Delta t}+G)} }$
\end{tabular}
\end{minipage}
\end{minipage} 
\medskip
\newline
The method diverges if $|\rho({P})| >1$.  $L_1,C,G$ and $L_2$ are fixed so the convergence of the method depends on $\Delta t$. For $\Delta t_0 \stackrel{def}{=}\frac{L_1G+\sqrt{(L_1G)^2+4L_1C}}{2} $, we have $|\rho({P})|=1$.
 The method converges with choosing a $\Delta t \in \rbrack 0 ; \Delta t_0 \lbrack$, stagnates if $ \Delta t=\Delta t_0 $ and diverges otherwise.

\begin{figure}[h]
\begin{minipage}{14cm}
\begin{minipage}{6.8cm}
\includegraphics[scale=0.280]{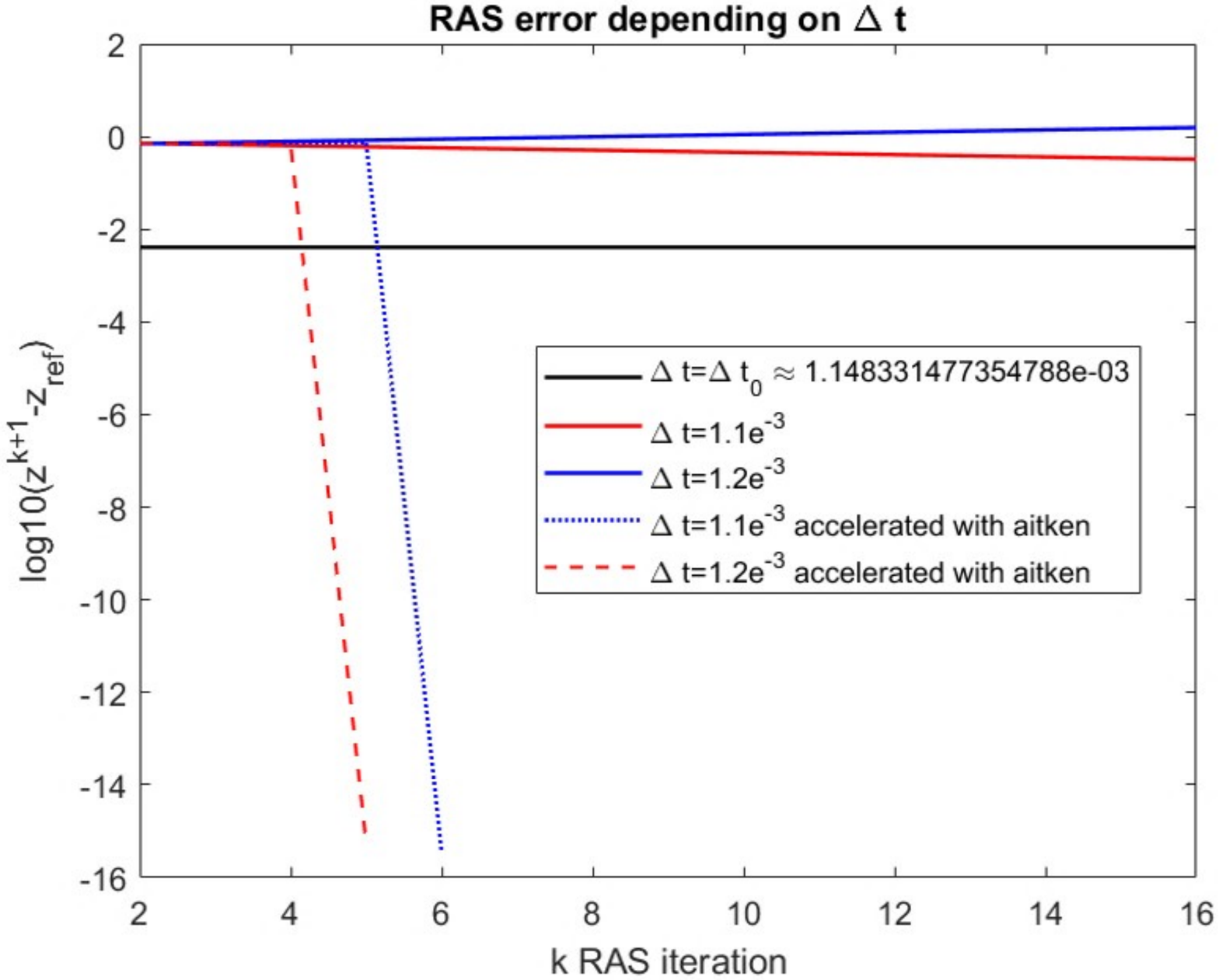}
\end{minipage}
\hfill
\begin{minipage}{6.8cm}
\includegraphics[scale=0.5]{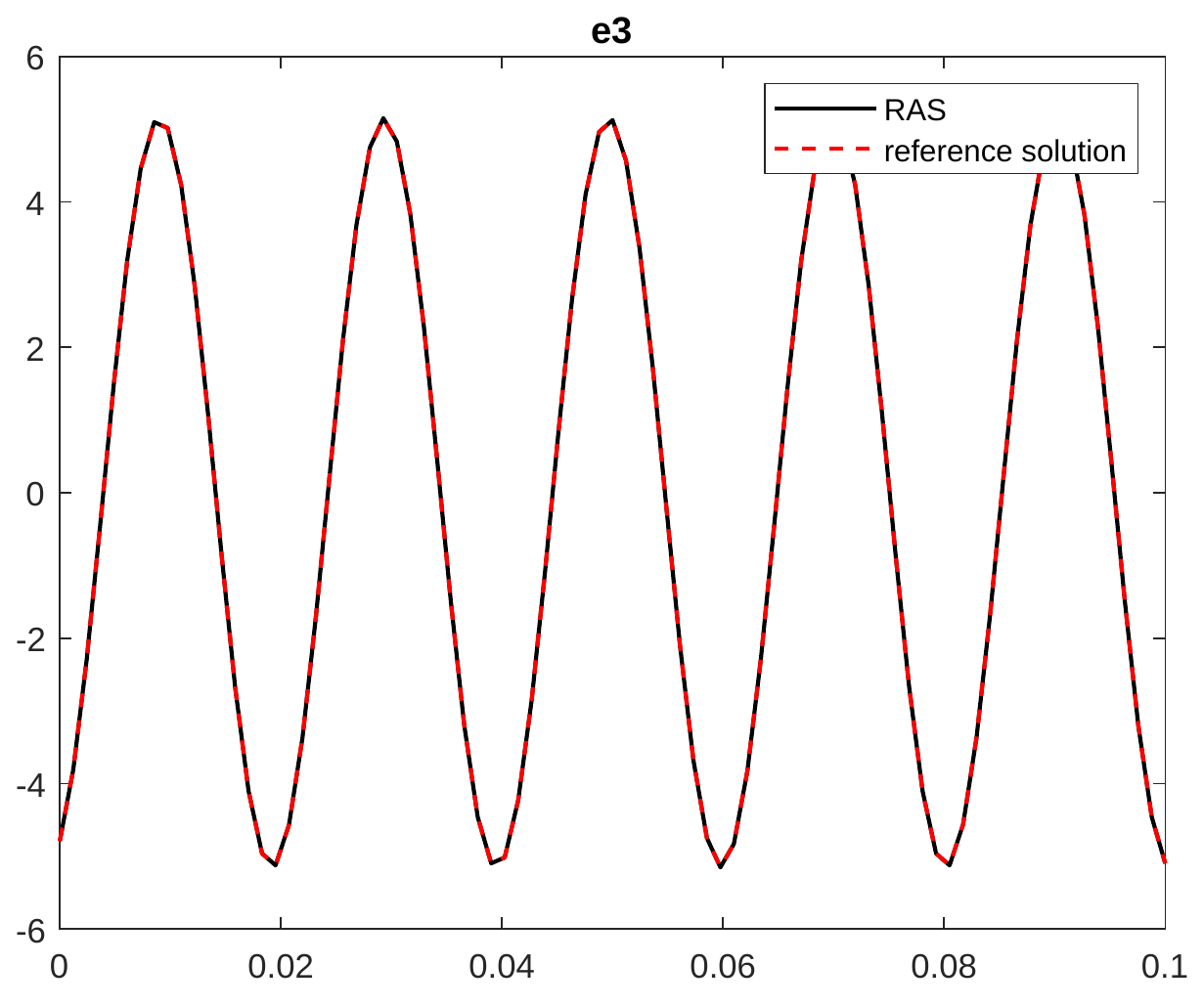}
\end{minipage}
\end{minipage}
\caption{First example with $ L_1=0.4, \,L_2=0.5, \,C=1.10^{-6}, \,G=2.10^{-3}$: (Left) DI with the RAS splitting  convergence behavior ( $log_{10}( ||z^{(2k)}-z_{ref}||_\infty)$  ) with respect to the RAS iterations and (right) comparison between the DI with the RAS splitting with the Aitken's technique for accelerating convergence and the DAE monolithic reference for the $e_3$ variable with $\Delta t_=1.2.10^{-3}$. }\label{shourick_contrib_Fig1}

\end{figure}

 Figure \ref{shourick_contrib_Fig1} (Left) gives the convergence behavior of  $log_{10}( ||z^{(2k)}-z_{ref}||_\infty)$ with respect to the RAS iterations for one time step for the three time step value cases while Figure \ref{shourick_contrib_Fig1}  (right) gives the $e_3$ behavior with respect to the time with the monolithic reference and the DI with the RAS splitting with the   Aitken's acceleration. 
\subsection{Second example of \cite{Pade_Tischendorf_WR_CVGCriterion_NUMAlgo_2019}}
For this second  example, the criteria of  \cite{Pade_Tischendorf_WR_CVGCriterion_NUMAlgo_2019} cannot ensure the convergence for the WR Gauss Seidel method. The circuit splitting is as follows:
\begin{figure}[h]
\begin{minipage}{10.5cm}
\begin{minipage}{5.2cm}
 \begin{tikzpicture}[scale=0.55]
 \draw[Cerise] (-6,2) -- (-6,1.15); 
  \draw[DodgerBlue]   (-6,1) -- (-6,-2); 
 \draw[Cerise]  (-6,2) -- (-3,2);
  \draw[Cerise] (-2,2) -- (1,2);
 \draw[Cerise]  (1,2) -- (1,1.15); 
  \draw[DodgerBlue]   (1,1) -- (1,-2); 
 
  \draw[DodgerBlue]   (1,1)--(-0.5,1); 
    \draw[DodgerBlue] (-1.5,1) -- (-3.5,1); 
   \draw[DodgerBlue] (-4.5,1) -- (-6,1); 
   
    \draw[DodgerBlue] (-2.5,1) -- (-2.5,-0.3); 
     \draw[DodgerBlue] (-2.5,-0.7) -- (-2.5,-2); 
     
       \draw[DodgerBlue] (1,-2)--(-0.5,-2); 
    \draw[DodgerBlue] (-1.5,-2)-- (-3.5,-2); 
   \draw[DodgerBlue] (-4.5,-2) -- (-6,-2); 
     
  \fill[DodgerBlue] (-6,1) circle(0.07);
   \fill[DodgerBlue] (1,1) circle(0.07);
  \fill[DodgerBlue] (-2.5,-2) circle(0.07);
   \fill[DodgerBlue] (-2.5,1) circle(0.07);
   
  \draw[DodgerBlue] (-6,1) node[left]{\scriptsize $n_1$};
   \draw[DodgerBlue] (1,1) node[right]{\scriptsize $n_r$};
   \draw[DodgerBlue]  (-2.5,-2) node[below]{\scriptsize $n_3$};
   \draw[DodgerBlue] (-2.5,1) node[above]{\scriptsize $n_2$};
   
   \begin{scope}[shift={(-3.2,2)},rotate=90]
{
 \draw[white, fill=white] (-0.1,0) rectangle (0.4,-1.5);
 \foreach \r in {0,...,2}
 {
  \draw[Cerise,thick,scale=1/3,shift={(0,-\r)}]
	(0,0) .. controls ++(2,0) and ++(1,0) ..
	++(0,-1.5) .. controls ++(-1,0) and ++(-0.5,0) ..
	++(0,0.5);
 }
 \draw[Cerise,thick,scale=1/3,shift={(0,-3)}] (0,0) .. controls ++(2,0) and ++(1,0) .. ++(0,-1.5);
}

\end{scope}
\draw[Cerise]  (-2.3,2.4) node[above]{\footnotesize  $L_1$};

        \begin{scope}[shift={(-1.5,1)},rotate=90]
{
 \draw[white, fill=white] (-0.1,0) rectangle (0.4,-1.5);
 \foreach \r in {0,...,2}
 {
  \draw[DodgerBlue,thick,scale=1/3,shift={(0,-\r)}]
	(0,0) .. controls ++(2,0) and ++(1,0) ..
	++(0,-1.5) .. controls ++(-1,0) and ++(-0.5,0) ..
	++(0,0.5);
 }
 \draw[DodgerBlue,thick,scale=1/3,shift={(0,-3)}] (0,0) .. controls ++(2,0) and ++(1,0) .. ++(0,-1.5);
}
\end{scope}
 \draw[DodgerBlue] (-0.6,0.25) node[above]{\footnotesize $L_2$};

   \draw[DodgerBlue,ultra thick] (-2.1,-0.3) -- (-2.8,-0.3);
   \draw[DodgerBlue,ultra thick] (-2.1,-0.7) -- (-2.8,-0.7);

  \draw[DodgerBlue] (-1.9,-0.5) node{\footnotesize  $C$};

   \draw[DodgerBlue,thick] (-4,-2) circle(0.5);
  \draw[DodgerBlue,thick] (-3.5,-2)-- (-4.5,-2);

  \draw[DodgerBlue] (-4,-2.8) node{\footnotesize $E_v $};

   \draw[DodgerBlue,thick] (-1,-2) circle(0.5);
  \draw[DodgerBlue,thick] (-1,-2.5)-- (-1,-1.5);

  \draw[DodgerBlue] (-1,-2.8) node{\footnotesize $E_i $};

    \draw[DodgerBlue](-4.5,1.15) -- (-3.5,1.15); 
   \draw[DodgerBlue] (-4.5,0.85) -- (-3.5,0.85); 
    \draw[DodgerBlue](-4.5,1.15) -- (-4.5,0.85); 
   \draw[DodgerBlue] (-3.5,1.15) -- (-3.5,0.85);

   \draw[DodgerBlue](-4,0.25) node[above]{\footnotesize $G$};
 
\end{tikzpicture}
\vspace*{-0.5cm}
\begin{eqnarray}
\small
\color{DodgerBlue}e_r&\color{DodgerBlue}=&\color{DodgerBlue}0,\nonumber \\
\color{DodgerBlue}e_1-e_3-E_{v}(t)-Z_s i_{5} &\color{DodgerBlue}=&\color{DodgerBlue} 0,\nonumber \\
\color{Cerise}e_1-e_r-L_1 \dfrac{d{i}_{1}}{dt} &\color{Cerise}=& \color{Cerise}0,\nonumber 
\end{eqnarray}
\end{minipage}
\hfill
\begin{minipage}{5.4cm}
\begin{eqnarray}
\small
\color{DodgerBlue}G (e_1-e_2) -i_{2} &\color{DodgerBlue}=&\color{DodgerBlue} 0, \nonumber \\
\color{DodgerBlue}(e_2-e_r)-L_2 \dfrac{d{i}_{3}}{dt} &\color{DodgerBlue}=& \color{DodgerBlue}0, \nonumber \\
\color{DodgerBlue}C (\dfrac{d{e}_3}{dt}-\dfrac{d{e}_2}{dt})-i_{4} &\color{DodgerBlue}=&\color{DodgerBlue}0, \nonumber\\
\color{DodgerBlue}E_{i}(t)-i_{6} &\color{DodgerBlue}=& \color{DodgerBlue}0, \nonumber \\
\color{DodgerBlue}-i_{1}-i_{6}-i_{3}&\color{DodgerBlue}=&\color{DodgerBlue} 0, \nonumber \\
\color{DodgerBlue}i_3-i_2-i_4&\color{DodgerBlue}=& \color{DodgerBlue}0, \nonumber \\
\color{DodgerBlue}-i_{5}+i_{4}+i_{6}&\color{DodgerBlue}=&\color{DodgerBlue} 0. \nonumber 
\end{eqnarray}
\end{minipage}
\end{minipage}
\end{figure}
The interface values are  $e_1$ and $e_r$ for the first partition and $i_1$ for the second partition.
The error operator ${P}$ and it's eigen values are calculated following Eq. \eqref{contribCalcul0} and \eqref{contribCalcul1}:\\
\begin{minipage}{14cm}
\begin{minipage}{7cm}
\begin{eqnarray}
{P}&=&\left(\begin{array}{ccc}
0 &\colorboxed[rgb]{0.12, 0.56, 1.0}{\begin{matrix}  -\frac{\Delta t}{L_1} &\frac{\Delta t}{L_1} \end{matrix}}\ \\
\colorboxed[rgb]{1, 0.3, 1.0}{\frac{\frac{L_2C}{\Delta t^2}+\frac{L_2G}{\Delta t}+1}{\frac{C}{\Delta t}+G}}&0&0\\
0&0&0
\end{array}  \right),  \nonumber
\end{eqnarray}
\end{minipage}
\hfill
\begin{minipage}{6cm}
\begin{tabular}{c|c}
$\lambda_1$ & 0\\
\hline
$\lambda_2$ & $-i \sqrt{\frac{L_2}{L_1}+\frac{\Delta t}{L_1(\frac{C}{\Delta t}+G)}}$\\
\hline
$\lambda_3$ & $i \sqrt{\frac{L_2}{L_1}+\frac{\Delta t}{L_1(\frac{C}{\Delta t}+G)} }$
\end{tabular}
\end{minipage}
\end{minipage} 
 The same way as for the first example, the method diverges if $|\rho({P})| >1$. So if $L_2\ge L_1$ the methode diverges,  $L_1,C,G$ and $L_2$ are fixed the way as $L_2<L_1$, so the convergence of the method depends on $\Delta t$. For $\Delta t_0 \stackrel{def}{=} \frac{(L_1-L_2)G+\sqrt{(L2-L1)^2G^2+4(L_1-L_2)C}}{2}$, we have $|\rho({P})|=1$.
 The method converges with choosing a $\Delta t \in \rbrack 0 ; \Delta t_0 \lbrack$, stagnates if $ \Delta t=\Delta t_0 $ and diverges otherwise.
 \begin{figure}[h!]
\begin{minipage}{13.6cm}
\begin{minipage}{6.8cm}
\includegraphics[scale=0.25]{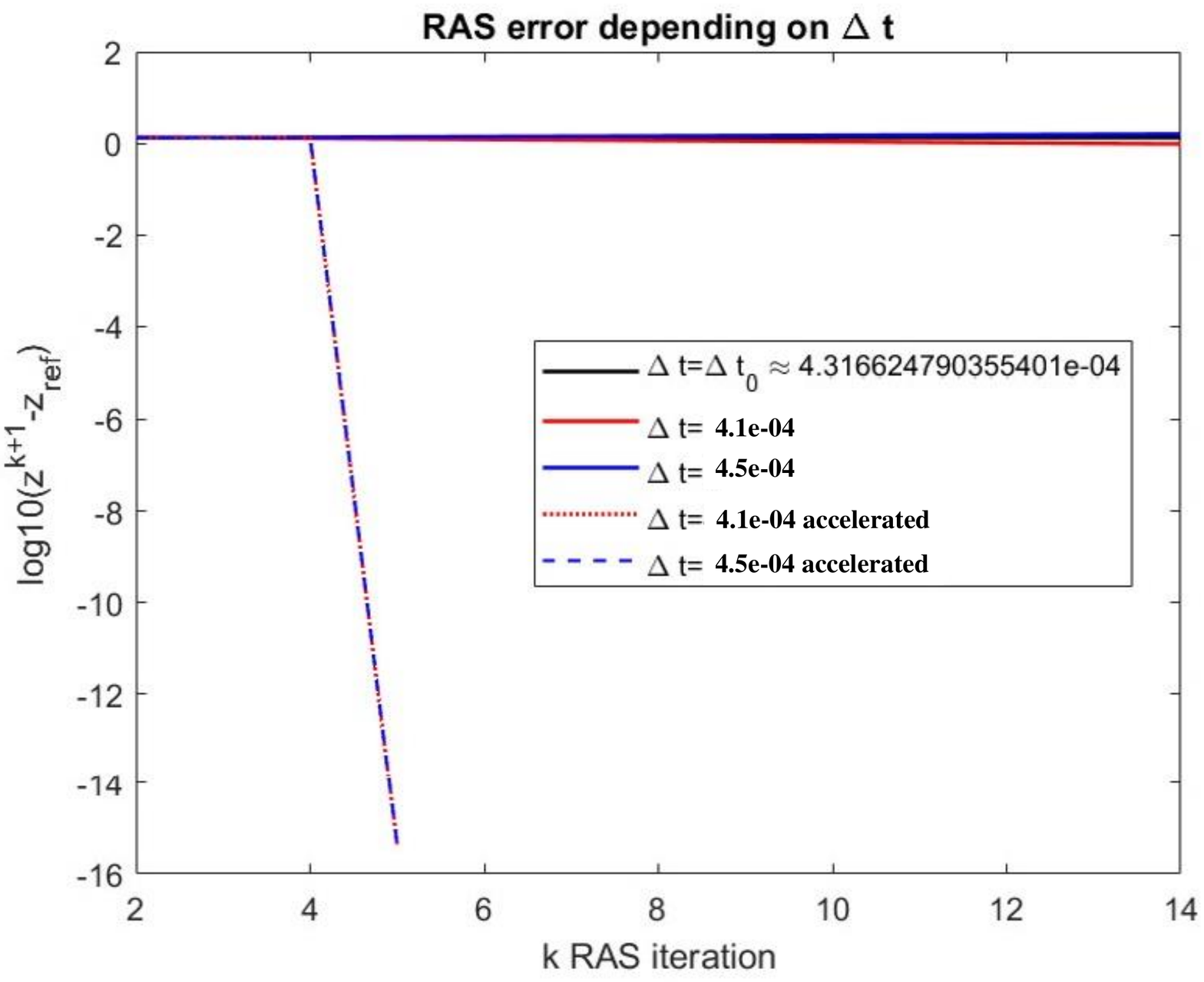}
\end{minipage}
\hfill
\begin{minipage}{6.8cm}
\includegraphics[scale=0.25]{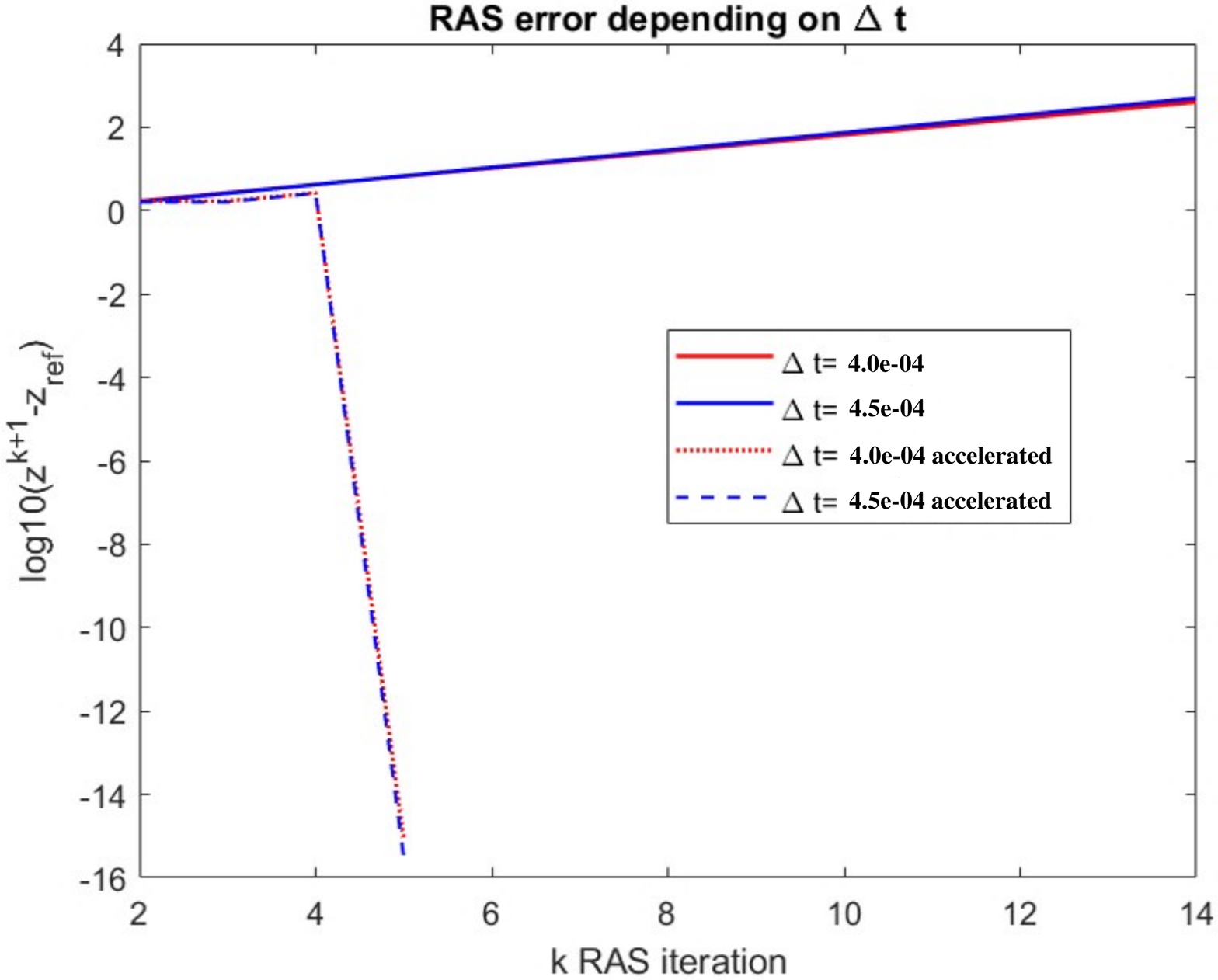}
\end{minipage}
\end{minipage}
\caption{ Second example:  DI with the RAS splitting convergence behavior ( $log_{10}( ||z^{2k}-z_{ref}||_\infty)$) with respect to the RAS iterations and  the Aitken's technique for accelerating convergence   after four RAS iterates, plus one more local solving,  with, $ C=1.10^{-6},\,G=2.10^{-3}$   and  (right) $ L_1=0.4 > L_2=0.3$,  (Left) $ L_1=0.5 < L_2=0.7$}\label{shourick_contrib_Fig2}
\end{figure}

 Figure \ref{shourick_contrib_Fig2} (Left) gives the convergence behavior of $log_{10}( ||z^{(2k)}-z_{ref}||_\infty)$ with respect to the RAS iterations for one time step for the three time step value cases and  the Aitken's technique for accelerating convergence   after four RAS iterates plus one more local solving, with $L_2< L_1$. It  shows that in both cases convergent or divergent the Aitken's acceleration reaches the monolithic reference  solution.  Figure \ref{shourick_contrib_Fig2} (right) gives the convergence behavior of  $log_{10}( ||z^{(2k)}-z_{ref}||_\infty)$ with respect to the RAS iterations  and  the Aitken's technique for accelerating convergence  after four RAS iterates plus one more local solving, with $L_2> L_1$. In the two time step cases the RAS diverges, but the Aitken's acceleration succeeds to reach the monolithic reference solution.
 
\begin{figure}[h!]
 \centering
\includegraphics[scale=0.5]{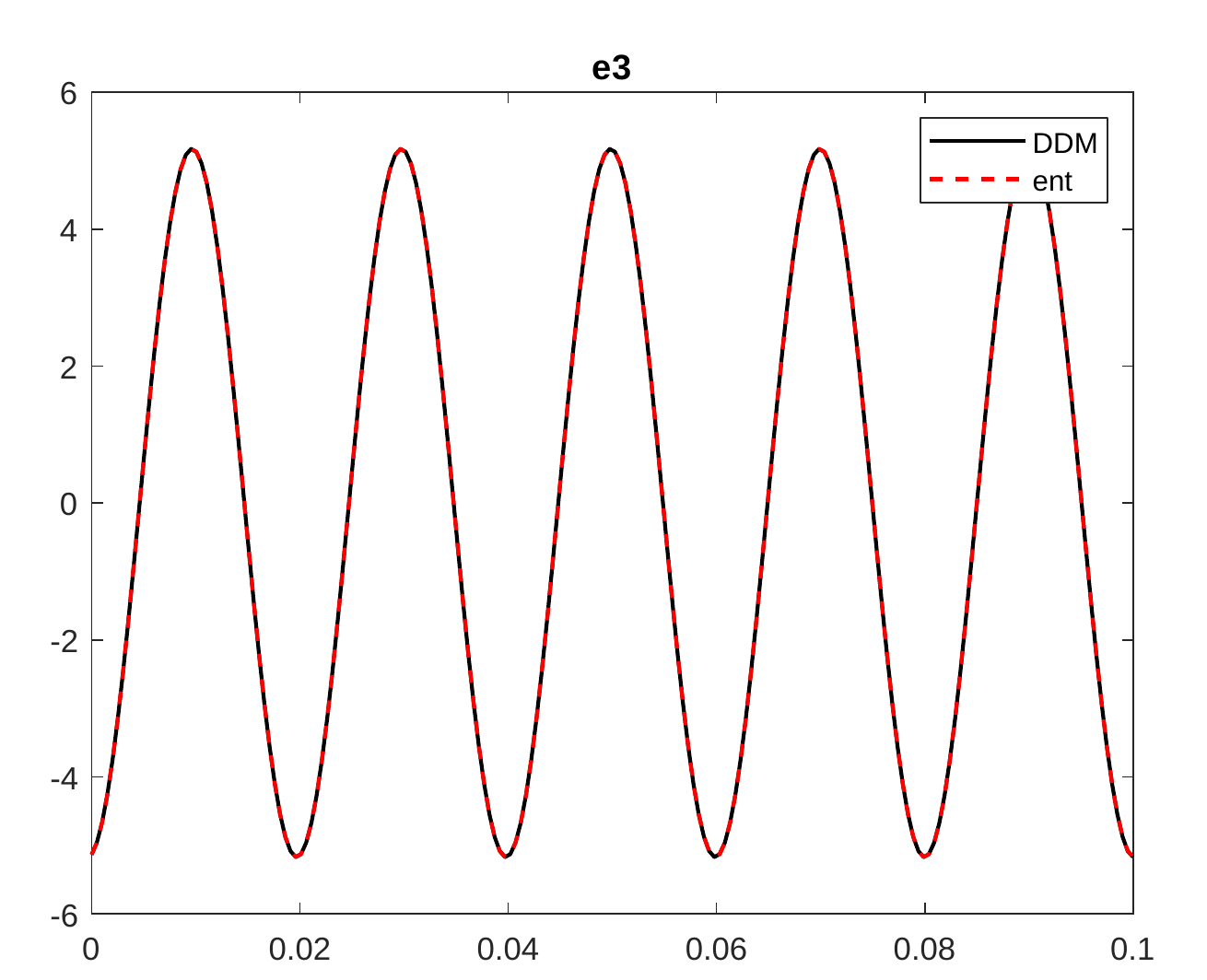}
\caption{ Second example:  Comparison between the  DI with the RAS splitting with the Aitken's technique for accelerating convergence and the DAE monolithic reference for the $e_3$ variable with $\Delta t =4.5 10^{-4}$ and for  $ C=1.10^{-6},\,G=2.10^{-3},\,  L_1=0.5 < L_2=0.7$. }\label{shourick_contrib_Fig2b}
\end{figure}

 Figure \ref{shourick_contrib_Fig2b}  gives the $e_3$ behavior with respect to the time with the monolithic reference solution and the DI with the RAS splitting and with the Aitken's acceleration, for  $ C=1.10^{-6},\,G=2.10^{-3},\,  L_1=0.5 < L_2=0.7$.

\subsection{DI with RAS splitting with heterogeneous modeling}
Secondly, we apply the method to the co-simulation of an RLC circuit split into two different types of modeling (ElectroMagnetic Transient (EMT): a very precise model requiring very small time steps and a dynamic phasor (TS) model: less precise but allowing the use of larger time steps). 
 The numerical example is the   RLC circuit of Figure \ref{shourick_contrib_Fig3ddm}.
  \begin{figure}[h!]
\begin{minipage}{13.5cm}
\begin{minipage}{7.2cm}
 \begin{tikzpicture}[scale=0.53]
  \draw[black!60!DodgerBlue!70] (-3.65,2.7) node[above]{$W^1_0$};
 \draw[black!70!yellow!65!red!65!] (2.55,2.7) node[above]{$W^1_1$};

 \draw (-4,2) -- (-3,2);
 \draw (-2,2) -- (-1.4,2);
 \draw(0,2) -- (1.5,2);
 \draw (2,2) -- (3,2);
 \draw (3,2) -- (3,-2); 
 \draw (-4,-2) -- (-2.5,-2);
 \draw (-2,-2) -- (-1,-2);
 \draw (0,-2) -- (1,-2);
 \draw (2.4,-2) -- (3,-2);
 
  \draw (-4,2) -- (-4,-2);

   \fill(-4,2) circle(0.05);
   \fill(-1.6,2) circle(0.05);
   \fill (0.6,2) circle(0.05);
   \fill (3,2) circle(0.05);
   \fill (-1.6,-2) circle(0.05);
   \fill (0.6,-2) circle(0.05);
   \fill (-4,-2) circle(0.05);

  \draw (-4,2) node[above]{\scriptsize $2$};
   \draw(-1.6,2) node[above]{\scriptsize $3$};
   \draw (0.6,2) node[above]{\scriptsize $4$};
   \draw (3,2) node[above]{\scriptsize $5$};
   \draw (-1.6,-2) node[below]{\scriptsize $7$};
   \draw (0.6,-2) node[below]{\scriptsize $6$};
   \draw (-4,-2) node[left]{\scriptsize $1$};
  
 \draw[thick] (1.5,2.5) -- (1.5,1.5);
 \draw[thick] (2,2.5) -- (2,1.5);

  \draw (1.7,2.4) node[above]{\footnotesize  $C_1$};
   
 \draw[thick] (-2,-2.5) -- (-2,-1.5);
 \draw[thick] (-2.5,-2.5) -- (-2.5,-1.5);
 
  \draw (-2.3,-1.6) node[above]{\footnotesize $C_2$};

  \draw[ thick] (-4,-2) -- (-4,-2.5);
 \draw[thick] (-4.5,-2.5) -- (-3.6,-2.5);
  \draw[thick] (-4.5,-2.5) -- (-4.3,-2.7);
 \draw[thick] (-4.2,-2.5) -- (-4,-2.7);
  \draw[thick] (-3.9,-2.5) -- (-3.7,-2.7);
   \draw[thick] (-3.6,-2.5) -- (-3.4,-2.7);
   
   
    \draw[thick] (-1.4,2) -- (-1.2,2.3);
    \draw[thick] (-1.2,2.3) -- (-1,1.7);
     \draw[thick] (-1,1.7) -- (-0.8,2.3);
 \draw[thick] (-0.8,2.3) -- (-0.6,1.7);
  \draw[thick] (-0.6,1.7) -- (-0.4,2.3);
  \draw[thick] (-0.4,2.3) -- (-0.2,1.7);
  \draw[thick] (-0.2,1.7) -- (0,2);
  
   \draw (-0.7,2.4) node[above]{\footnotesize $R_1$};
  
  \draw[thick] (1,-2) -- (1.2,-1.7);
  \draw[thick] (1.2,-1.7) -- (1.4,-2.3);
 \draw[thick] (1.4,-2.3) -- (1.6,-1.7);
 \draw[thick] (1.6,-1.7) -- (1.8,-2.3);
  \draw[thick] (1.8,-2.3) -- (2,-1.7);
  \draw[thick] (2,-1.7) -- (2.2,-2.3);
  \draw[thick] (2.2,-2.3) -- (2.4,-2);
   \draw (1.7,-1.6) node[above]{\footnotesize $R_2$};
  generateur de tension
  
  \draw[thick] (-4,0) circle(0.5);
  \draw[thick] (-4,0.5)-- (-4,-0.5);
  \draw[->,thick] (-3.45,-0.5)-- (-3.45,0.5);
  \draw (-3.45,0) node[right]{\footnotesize E cos $\omega t = \beta$};

  inductance
   
   \begin{scope}[shift={(-3.5,2)},rotate=90]
{
 \draw[black!7!DodgerBlue!11,opacity=0.80,fill=black!7!DodgerBlue!11,opacity=0.80] (-0.2,0) rectangle (0.4,-1.5);
 3 spires definies par des courbes de bezier, via une boucle for
 \foreach \r in {0,...,2}
 {
  \draw[thick,scale=1/3,shift={(0,-\r)}]
	(0,0) .. controls ++(2,0) and ++(1,0) ..
	++(0,-1.5) .. controls ++(-1,0) and ++(-0.5,0) ..
	++(0,0.5);
 }
 une demi-spire pour finir
 \draw[thick,scale=1/3,shift={(0,-3)}] (0,0) .. controls ++(2,0) and ++(1,0) .. ++(0,-1.5);
}

\end{scope}
\draw (-2.8,2.4) node[above]{\footnotesize  $L_1$};

   \begin{scope}[shift={(-1.2,-2)},rotate=90]
{
 \draw[black!5!yellow!10!red!8!,opacity=0.3,fill=black!5!yellow!10!red!8!,opacity=0.8] (-0.2,0) rectangle (0.4,-1.5);
  \draw[black!5!yellow!10!red!8!,opacity=0.3,fill=black!5!yellow!10!red!8!,opacity=0.7] (-0.2,0) rectangle (0.25,-1.2);
 3 spires definies par des courbes de bezier, via une boucle for
 \foreach \r in {0,...,2}
 {
  \draw[thick,scale=1/3,shift={(0,-\r)}]
	(0,0) .. controls ++(2,0) and ++(1,0) ..
	++(0,-1.5) .. controls ++(-1,0) and ++(-0.5,0) ..
	++(0,0.5);
 }
 une demi-spire pour finir
 \draw[thick,scale=1/3,shift={(0,-3)}] (0,0) .. controls ++(2,0) and ++(1,0) .. ++(0,-1.5);
}
\end{scope}
 \draw (-0.5,-1.6) node[above]{\footnotesize $L_2$};

\draw[blue,thick,->] (-4.5,2.6) -- (-4,2.6);
\draw[red,thick,<-] (0.9,-1.3) -- (1.3,-1.3);
\draw[blue](-4.4,2.45) node[above]{\footnotesize ${\bf i_{12},v_2}$};
\draw[red](1.3,-1.1) node[below]{\footnotesize ${\bf i_{56},v_6}$};

\coordinate (h) at (-3.8,2.5);
\coordinate (i) at (3.5,2.5);
\coordinate (j) at (3.5,-2.3);
\coordinate (k) at (-1.9,-2.3);

\coordinate (m) at (-4.4,2.2);
\coordinate (mn) at (-1.5,2.8);
\coordinate (n) at (0.63,2.2);
\coordinate (no) at (1,0);
\coordinate (o) at (0.32,-2.2);
\coordinate (op) at (-2,-2.9);
\coordinate (p) at (-4.4,-2.2);
\coordinate (pm) at (-4.85,0);

  \draw[dotted,draw=red!60!yellow!35!black,fill=black!5!yellow!10!red!15!,opacity=0.50](h) .. controls +(1,0.8) and +(-0.9,0.8) .. (i)
               .. controls +(0.8,-0.9) and +(0.7,0.9) .. (j)
               .. controls +(-0.9,-0.8) and +(0.9,-0.7) .. (k)
               .. controls +(-1,1) and +(0.9,-1) .. (h); 

\draw[dotted,draw=black,fill=black!5!DodgerBlue!11,opacity=0.5]  (m) .. controls +(1,0.5) and +(-1,0.08) .. (mn)
               .. controls +(1,0.08) and +(-0.8,0.5) .. (n)
               .. controls +(0.5,-0.8) and +(-0.17,2) .. (no)
               .. controls +(-0.17,-2) and +(0.5,0.8) .. (o)
               .. controls +(-0.7,-0.5) and +(1,0.1) .. (op)
               .. controls +(-1,0.1) and +(0.7,-0.5) .. (p)
                .. controls +(-0.5,0.8) and +(0.1,-0.8) .. (pm)
               .. controls +(0.1,0.8) and +(-0.5,-0.8) .. (m); 
\end{tikzpicture}
\vspace*{-0.5cm}
\begin{eqnarray}
\small
v_1&=&0,\nonumber \\
v_2-v_1-E-Z_s i_{12} &=& 0,\nonumber  \\
v_3-v_2-L_1 \dfrac{d{i}_{23}}{dt} &=& 0, \nonumber \\
v_4-v_3 -R_1 i_{34} &=& 0, \nonumber \\
i_{67}-i_{71}&=& 0. \nonumber 
\end{eqnarray}
\end{minipage}
\hfill
\begin{minipage}{6.4cm}
\begin{eqnarray}
\small
C_1 (\dfrac{d{v}_5}{dt}-\dfrac{d{v}_4}{dt})-i_{45} &=&0,\nonumber \\
v_6-v_5-R_2 i_{56} &=& 0,\nonumber \\
v_7-v_6-L_2 \dfrac{d{i}_{67}}{dt} &=& 0,\nonumber \\
C_2 (\dfrac{d{v}_1}{dt}-\dfrac{d{v}_7}{dt})-i_{71} &=& 0, \nonumber  \\
i_{12}-i_{23}&=& 0, \nonumber  \\
i_{23}-i_{34}&=& 0, \nonumber  \\
i_{34}-i_{45}&=& 0, \nonumber  \\
i_{45}-i_{56}&=& 0, \nonumber  \\
i_{56}-i_{67}&=& 0, \nonumber  
\end{eqnarray}
\end{minipage}
\end{minipage}
\caption{Linear RLC circuit and its associated  EMT modeling  DAE system with $W^0_0=\left\{ v_2, i_{12},  v_3, i_{34}, i_{45}, i_{56}, v_6, i_{71}\right\}$ and $W^0_1=\left\{ v_1, i_{23}, v4, v5, i_{67}, v_7\right\}$ . $L_1=L_2=0.7$, $C_1=C_2=1.10^{-6}$, $R_1=R_2=77$, $Z_s=1.10^{-6}$, $\omega=2\pi\,50$, $E=5$. }\label{shourick_contrib_Fig3ddm}
\end{figure}
  
An overlap is defined and the EMT equations from $W^1_1$ are changed in TS equations and solved for the $-1$, $1$ and $0$ dynamic phasor modes. The values to be exchanged are $i_{12},v_{2}$ from the EMT to the TS side and $i_{56},v_{6}$ from the TS side to the EMT one.

The two main difficulties to carry out the co-simulation reside in the difference in the representations of the variables and the time step difference. We choose $\Delta t_{ts}=m \Delta t_{emt} \; m\in \mathbb{N}$ and the interface values are exchanged at each TS time step. For TS modeling, the variables are assumed to oscillate with a specific angular frequency $\omega_0= \dfrac{2 \pi}{T}$ (where $T$ is the period) and its selected harmonics (dynamic phasor modes) taken from a subset $I=\left\{\ldots,-1,0,1, \ldots\right\}$:
\begin{equation}
z(t)= \sum_{k \in I} z_k(t) exp(i k \omega_0 t),\, z=\left\{x,y\right\}.
 \label{shourick_contrib_Eq3}
\end{equation}
 Introducing \eqref{shourick_contrib_Eq3} into \eqref{EqLinearDAE} leads after simplification (i.e orthogonality of the functions $exp(ik\omega_0 t)$  with respect to the dot product $[f,g]=\dfrac{1}{T} \int_{t}^{t+T} f(z)g(z)dz$) to another  DAE system that  takes into account the differential property of the dynamic phasors. The number of TS variables is then multiplied by the number of harmonics chosen, and the number of equations must be multiplied accordingly.
 
 Let's take back the \eqref{EqLinearDAERASsplitDiscret2} equation, it must be adapted to the subdomain solved with the TS or with the EMT modeling. First the subdomain solved with TS for the $n+1$ time step and the $k+1$ RAS iteration:
 
  \begin{eqnarray}
  \underbrace{\left(\begin{array}{c}  x_{ts}^{n+1,(k+1)}\\ y_{ts}^{n+1,(k+1)}\end{array} \right)}_{z_{ts}^{n+1,(k+1)}} &=& \underbrace{\left( \begin{array}{cc} \tilde{A}_{ts} & \tilde{B}_{ts} \\  C_{ts} & D_{ts} \end{array} \right)^{-1}}_{\tilde{\mathbb{A}}_{ts}^{-1}}   \underbrace{ \left( \left(\begin{array}{c} \tilde{b}_{ts,d}^{n+1} \\  b_{ts,a}^{n+1}
\end{array} \right) \right.}_{\tilde{b}^{n+1}_{ts}} - \underbrace{\left( \begin{array}{cc} \tilde{E}_{ts,d}^{d} & \tilde{E}_{ts,d}^{a} \\ E_{ts,a}^{d} &  E_{ts,a}^{a} \end{array} \right)}_{\tilde{\mathbb{E}}_{ts}} F_{\text{mod}} \underbrace{\left. \left( \begin{array}{c} X_{emt,e}^{n+1,(k)}  \\ Y_{emt,e}^{n+1,(k)} \end{array} \right) \right)}_ {Z_{emt,e}^{n+1,(k)}} 
 \nonumber
\end{eqnarray}
$F_{\text{mod}}$ represents a readjusted FFT and a choice of mods corresponding to those retained for the TS simulation. $ Z_{emt,e}^{n+1,(k)}$ is a history of values computed by the EMT subsystem during the previous RAS iteration (completed by some of the last values from the previous time steps if $T>\Delta t_{ts}$). This history is the size of a period and ends at the instant corresponding to the $n+1$ time step. $Z_{emt}^{n+1,(k)}=[z_{emt}^{n+1-j,(k)}, \hdots , z_{emt}^{n,(k)},z_{emt}^{n+1,(k)}]$\\

Secondly we perform the simulation for $W_0^1$ solved with EMT for each intermediate time step, the equation \eqref{EqLinearDAERASsplitDiscret2} is adapted to the subdomain solved, it gives for an intermediate time step $m+1$:
 \begin{eqnarray}
  \underbrace{\left(\begin{array}{c}  x_{emt}^{m+1,(k+1)}\\ y_{emt}^{m+1,(k+1)}\end{array} \right)}_{z_{emt}^{m+1,(k+1)}} &=& \underbrace{\left( \begin{array}{cc} \tilde{A}_{emt} & \tilde{B}_{emt} \\  C_{emt} & D_{emt} \end{array} \right)^{-1}}_{\tilde{\mathbb{A}}_{emt}^{-1}}   \underbrace{ \left( \left(\begin{array}{c} \tilde{b}_{emt,d}^{m+1} \\  b_{emt,a}^{m+1}
\end{array} \right) \right.}_{\tilde{b}^{m+1}_{emt}} - \underbrace{\left( \begin{array}{cc} \tilde{E}_{emt,d}^{d} & \tilde{E}_{emt,d}^{a} \\ E_{emt,a}^{d} &  E_{emt,a}^{a} \end{array} \right)}_{\tilde{\mathbb{E}}_{emt}} R_{\text{mod}}^{(m+1)} \underbrace{\left. \left( \begin{array}{c} x_{ts,e}^{n+1,(k)} \\ y_{ts,e}^{n+1,(k)} \end{array} \right) \right)}_ {z_{ts,e}^{n+1,(k)}}  
 \nonumber
 \end{eqnarray}
  \begin{figure}[t]
  \centering
  \begin{minipage}{13cm}
  \begin{minipage}{6.4cm}
\includegraphics[scale=0.45]{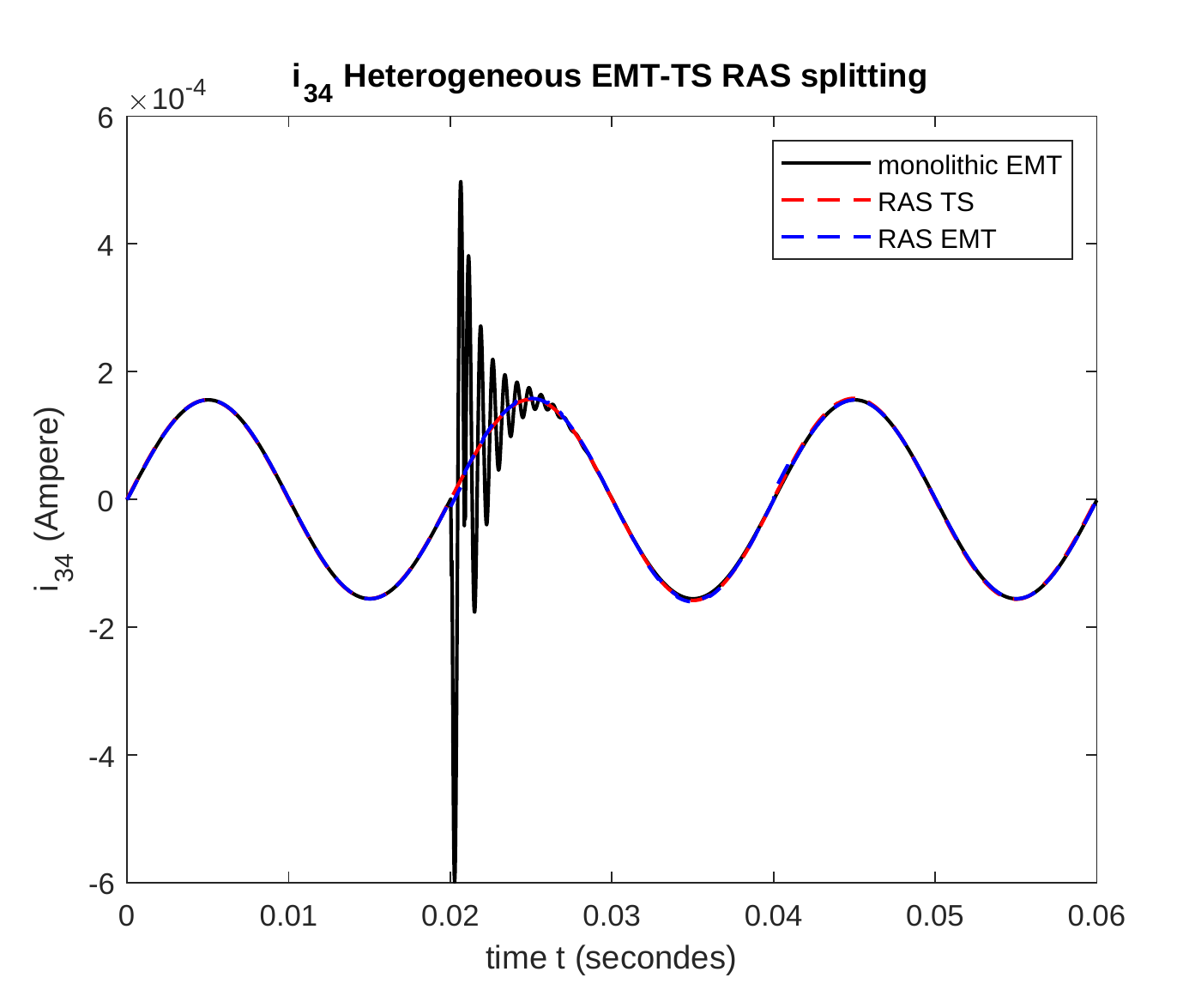}
\end{minipage}
\hfill
  \begin{minipage}{6.4cm}
\includegraphics[scale=0.45]{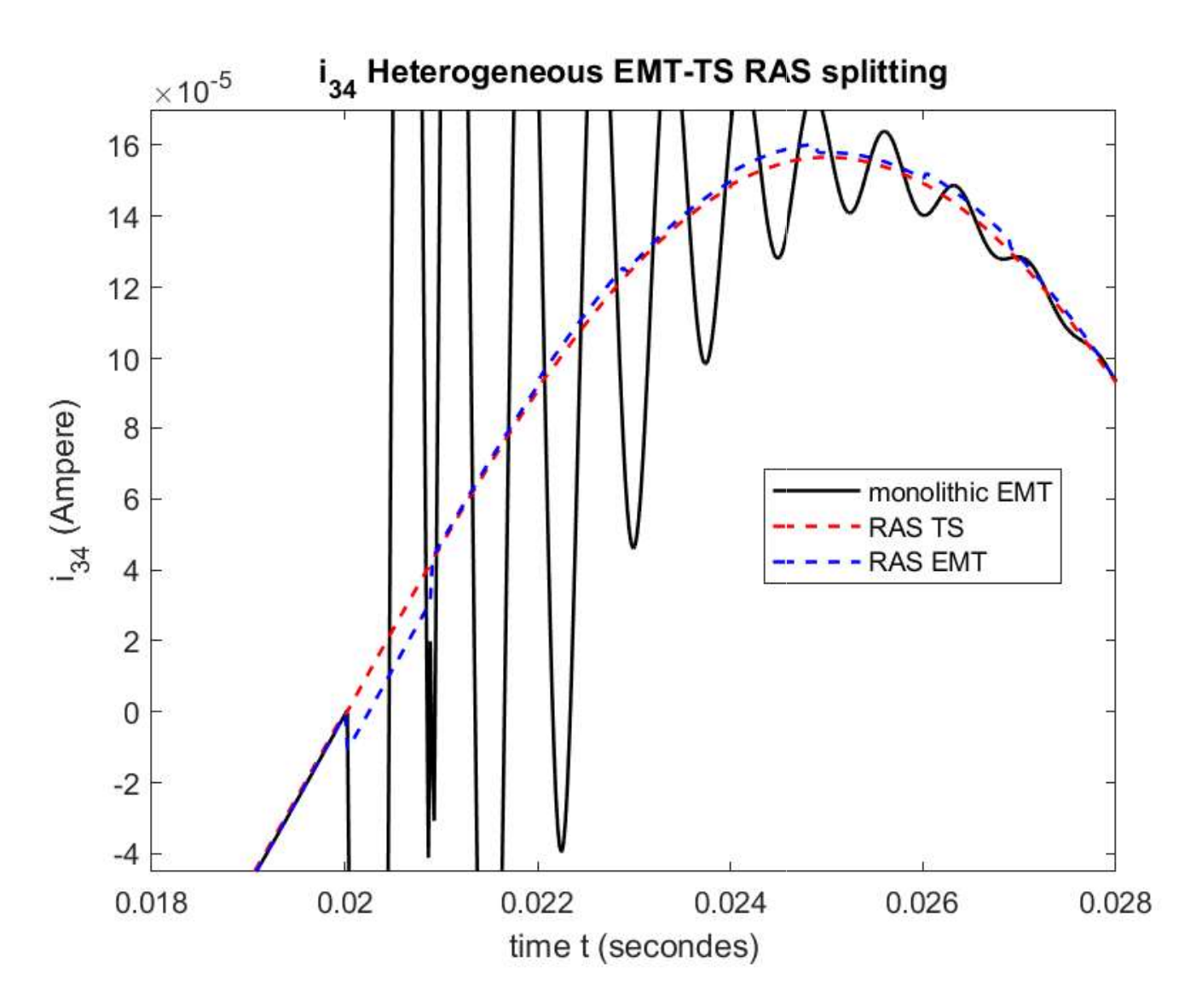}
\end{minipage}
\end{minipage}
   \begin{minipage}{13cm}
  \begin{minipage}{6.4cm}
\includegraphics[scale=0.45]{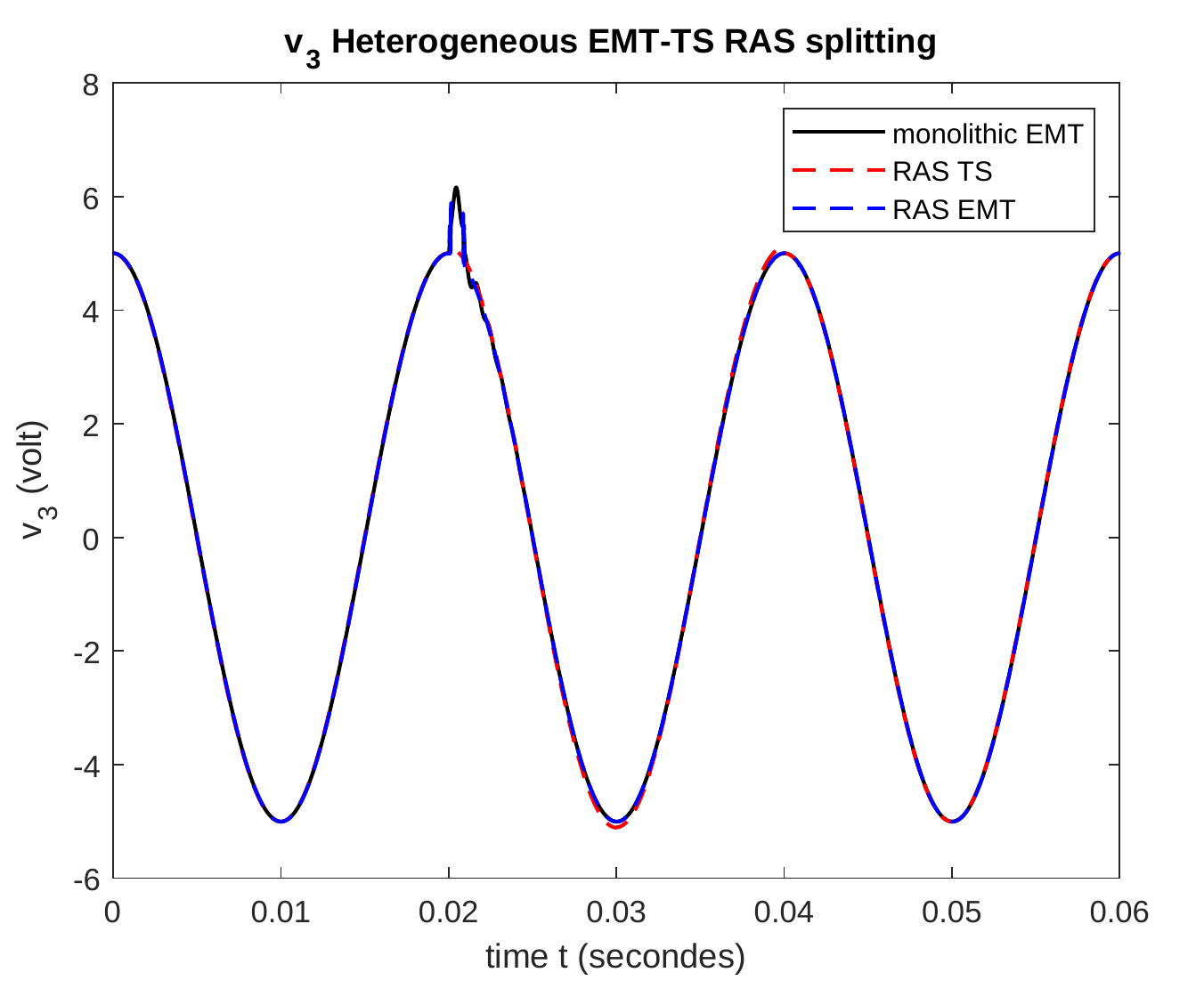}
\end{minipage}
\hfill
  \begin{minipage}{6.4cm}
\includegraphics[scale=0.45]{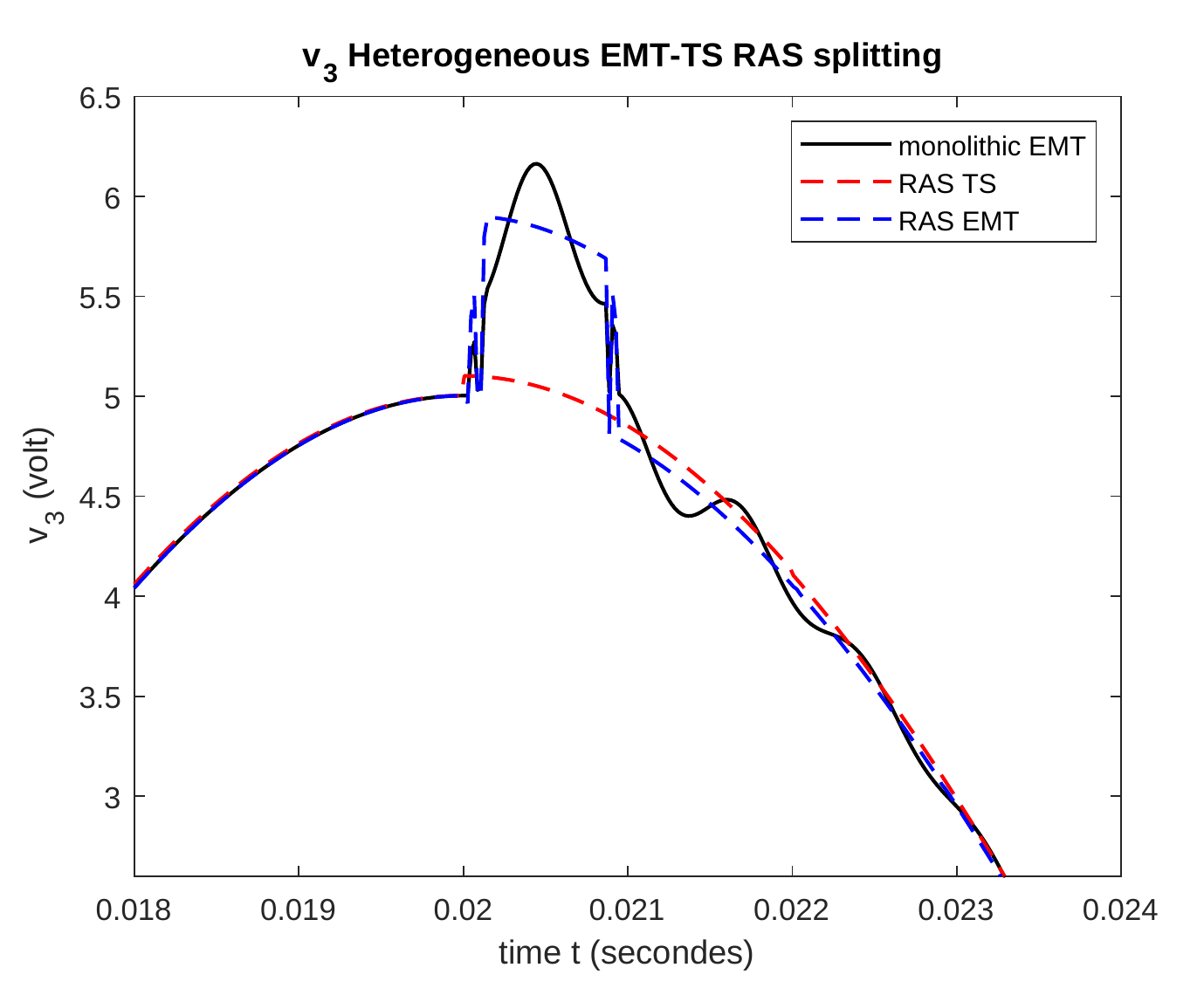}
\end{minipage}
\end{minipage}
\caption{ Comparison of the behavior, with respect to time, of the variables $i_{34}$ (top left) and $v_3$ (bottom left) (the figures on the right are their zoom on the disturbances) computed using the heterogeneous EMT-TS  RAS splitting with the Aitken's technique for accelerating convergence ($\Delta t_{ts}=2.10^{-3}$ and $\Delta t_{emt}=2.10^{-5}$), the reference is the monolithic EMT.  An amplitude perturbation on the voltage source starting  at $t=0.02$s and ending at $t=0.021$s, therefor lasting less than one $\Delta t_{ts}$ is applied. Parameters are $L_1=0.07, C_1=1.e-5, R_1=7, L_2=0.07, C_2=1.e-7, R_2=7, Zs=0.000001$.}\label{shourick_contrib_Fig3}
\end{figure}
 $z_{ts,e}^{n+1,(k)}$ are the values computed by the TS side at the $n+1$ time step, and $R_{\text{mod}}^{(m+1)}$ an operator which recombine the TS modes and estimate their values for the $m+1$ time step.
 \newline
  $F_{\text{mod}}$ and $R_{\text{mod}}$ are linear operators so the DI with the RAS splitting convergence/divergence always remains purely linear and so we can apply the Aitken's technique for accelerating convergence. Since the history $Z_{emt}^{n+1,(k)}$ can be very large, the resulting error matrix would be very cumbersome to invert. Therefore, the acceleration is only performed on the interface values computed by the TS side, then the converged TS interface values are used to resolve the EMT side locally and after the TS side locally.

Figure \ref{shourick_contrib_Fig3} compares the EMT monolithic reference values for the variables $v_3$ and $i_{34}$ with the EMT-TS heterogeneous RAS splitting where a perturbation on the source voltage that starts at $t=0.02$s and ends at $t=0.021$s is applied. The RAS DDM succeeds in capturing part of the perturbation on the $v_3$. It shows a good agreement between the monolithic and the EMT-TS heterogeneous RAS for the variable $v_3$ . The variable $i_{34}$ in the EMT DDM part captures certain oscillations due to the perturbation. These results show that EMT-TS heterogeneous RAS splitting can capture disturbances that last less than one TS time step and therefore would not have been captured by a monolithic TS model. 
\begin{figure}[h!]
\centering
\includegraphics[scale=0.6]{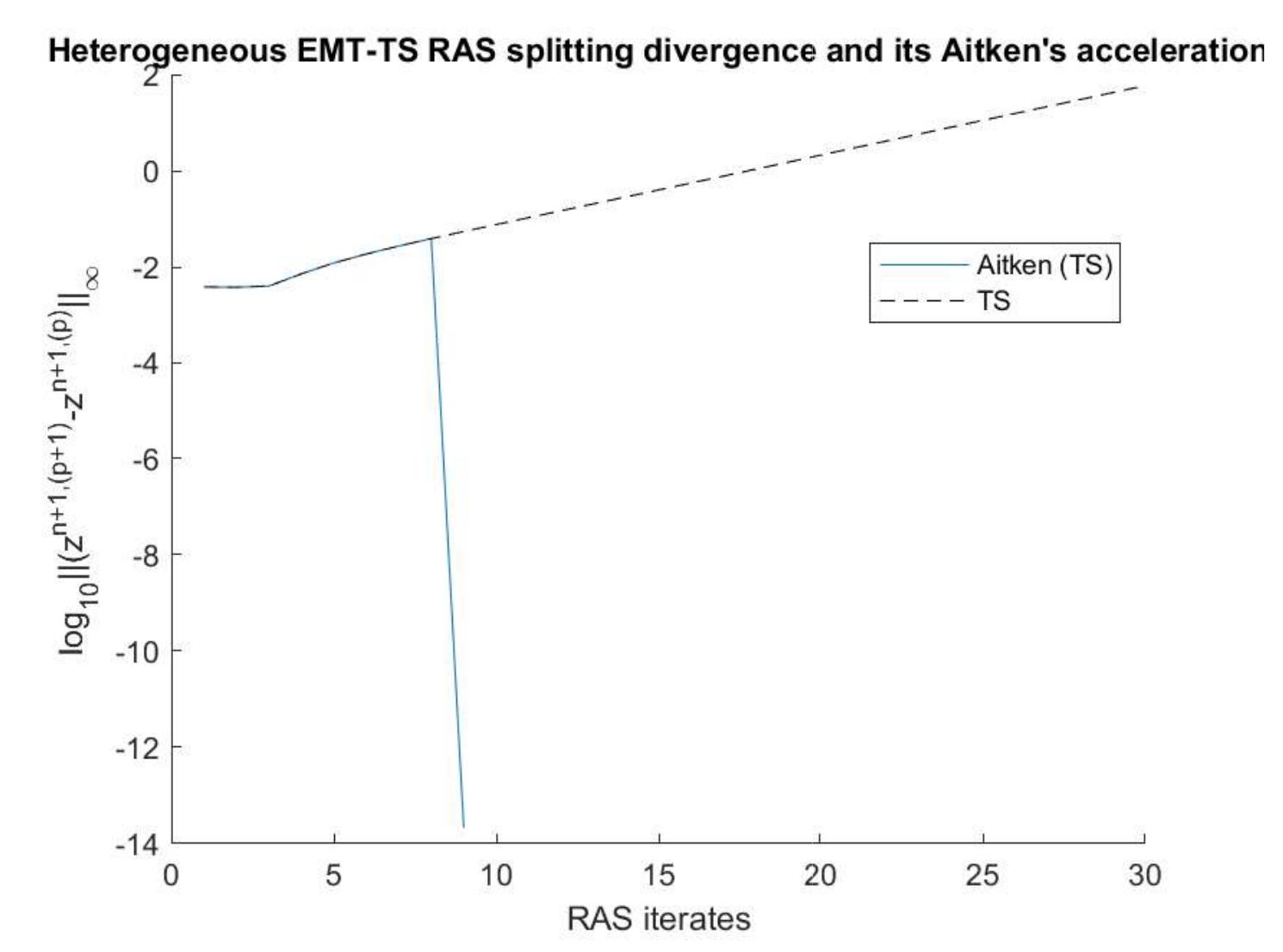}
\caption{(Left) Heterogeneous EMT-TS  RAS convergence error for the TS boundary for the time step $t = 0.02$ and its Aitken's acceleration applied on the TS partition interface  with $\Delta t_{ts}=2.10^{-3}$ and $\Delta t_{emt}=2.10^{-5}$) and with parameters $L_1=0.07, C1=1.e-6, R1=7, L2=0.07, C2=1.e-6, R2=7, Zs=0.000001$ }

\label{shourick_contrib_Fig3b}
\end{figure}

Figure \ref{shourick_contrib_Fig3b} shows the purely linear divergence of the DI with the RAS splitting and its  acceleration for the TS partition interface during the time step $t=0.02$  for  $\Delta t_{ts}=2.10^{-3}$ and $\Delta t_{emt}=2.10^{-5}$) with parameters $L_1=0.07, C1=1.e-6, R1=7, L2=0.07, C2=1.e-6, R2=7, Zs=0.000001$ chosen  to have DI with divergence.

\subsection{Non-linear Case}
We now consider, the problem \eqref{EqLinearDAE} but with at least one non-linear element. We rewrite it in its discrete form with linearizing it at each time step:
 \begin{eqnarray}
\underbrace{\left(\begin{array}{c}  x_i^{n+1,(k+1)}\\ y_i^{n+1,(k+1)}\end{array} \right)}_{z_i^{n+1,(k+1)}} &=& \underbrace{\left( \begin{array}{cc} \tilde{A}_i ^{n+1}& \tilde{B}_i^{n+1} \\  C_i^{n+1} & D_i^{n+1} \end{array} \right)^{-1}}_{(\tilde{\mathbb{A}}_i^{n+1})^{-1}}   \underbrace{ \left( \left(\begin{array}{c} \tilde{b}_{i,d}^{n+1} \\  b_{i,a}^{n+1}
\end{array} \right) \right.}_{\tilde{b}^{n+1}_i} - \underbrace{\left( \begin{array}{cc} \tilde{E}_{i,d}^{n+1,d} & \tilde{E}_{i,d}^{n+1,a} \\ E_{i,a}^{n+1,d} &  E_{i,a}^{n+1,a} \end{array} \right)}_{\tilde{\mathbb{E}}_i^{n+1}} \underbrace{\left. \left( \begin{array}{c} x_{ie}^{n+1,(k)}  \\ y_{ie}^{n+1,(k)} \end{array} \right) \right)}_ {z_{i,e}^{n+1,k}}   \nonumber
\end{eqnarray}
following the same steps as in section \ref {Dynamic_iteration_DDM}, the error between two iterations can be rewritten as:
\begin{eqnarray}
\left(
\begin{array}{c}
e_{\scriptscriptstyle 1}\\

e_{\scriptscriptstyle 2}\\

\end{array}  \right) ^{n+1,k+1}=
\left( \begin{array}{cc}
0& P_{\scriptscriptstyle1 } \\
P_{\scriptscriptstyle2 }& 0
\end{array}  \right)^{n+1}
\left(
\begin{array}{c}

e_{\scriptscriptstyle 1}\\

 e_{\scriptscriptstyle 2}\\

\end{array}  \right) ^{n,k}
\nonumber 
\end{eqnarray}
The $P^{n+1}$ error operator depends on the time step and so needs to be computed again for each time step. However, the error operator does not depend on the RAS iteration $k$ and so for each time step the (convergence/divergence) is purely linear and can be accelerated toward the true solution using the Aitken's technique for accelerating convergence.

Let us take again the numerical example 2 and replace G by a function of $ i_2 $ the current which crosses the associated component, $ G = \frac{1}{G_0 + \alpha i_2} $. Although this is non relevant from a physical point of view, we chose to take $ \alpha $ very large in order to increase the non-linearity. 

\begin{figure}[h]
\centering
\includegraphics[scale=0.36]{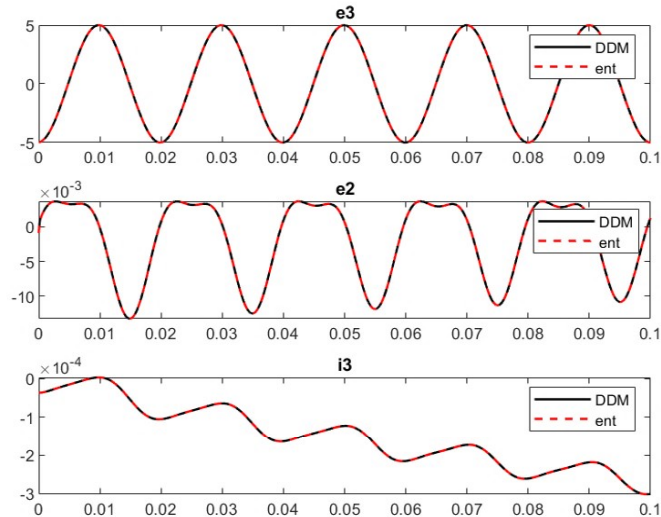}
\caption{Second numerical example with non linear G with $ L_1=0.6, \, L_2=0.7, C=1.10^{-6},\, G_0=10, \,\alpha=2000, \Delta t_=2.10^{-4}$: comparison between the DI with the RAS splitting with the Aitken's technique for accelerating convergence and the DAE monolithic reference for $e_2,e_3,i_3$ and (right) variations of the maximum eigen value of the error operator $P$ according to the time step.}\label{shourick_contrib_Fig5}
\end{figure}

\begin{table}[h]
\centering
\begin{tabular}{c|c|c|c}
\hline
$n$ &  1 & 25  & 250 \\
\hline
$\rho(P^n)$ & $1\pm 0816 i$ &$1\pm 0807i$ & $1\pm 0814 i$ \\
\hline
\end{tabular}
\caption{Variations of the maximum eigenvalue of the error operator $P^n$ according to the time step $n=\left\{ 1, 25, 250 \right\}$ for the second numerical example with non linear G with $ L_1=0.6, \, L_2=0.7, C=1.10^{-6},\, G_0=10, \,\alpha=2000, \Delta t=2.10^{-4}$. \label{shourick_contrib_Table}}
\end{table}
Table \ref{shourick_contrib_Table} gives the maximum eigenvalue of $ P^1$, $P^{25}$ and $P^{250}$.
It shows that the DI with the RAS splitting diverges for these time steps but with small variations in the maximum eigenvalue from one time step to another.\\
 Figure \ref{shourick_contrib_Fig5} shows identical behavior of  $e_2,e_3,i_3$, with respect to time, between the monolithic reference (dashed red curve) and the DI with the RAS splitting with the Aitken's acceleration (black curve).
  
\section{Conclusion \label{conclusion}}

We formulated the dynamic iteration method with the restricted additive Schwarz splitting as an iterative process involving the interface unknowns coming from the partitioning of the differential algebraic system of equations. Its pure linear convergence or divergence in the context of linear DAE system, allows us to accelerate the convergence toward the true solution with the Aitken's technique for accelerating convergence. We numerically built the error operator associated with the interface from the RAS iterations, only once if we use fixed time step. We also showed that the method can be used with heterogenous partition in the modeling such as EMT and TS modeling. Some extension of the method to solve nonlinear problems can also be applied by considering the linearization of the problem for each time step. 
\bibliographystyle{elsarticle-num}
\bibliography{DynamicIteration.bib}

\begin{thebibliography}{10}
\expandafter\ifx\csname url\endcsname\relax
  \def\url#1{\texttt{#1}}\fi
\expandafter\ifx\csname urlprefix\endcsname\relax\def\urlprefix{URL }\fi
\expandafter\ifx\csname href\endcsname\relax
  \def\href#1#2{#2} \def\path#1{#1}\fi

\bibitem{lelarasmee1982}
E.~Lelarasmee, A.~Ruehli, A.~Vincentelli, {The Waveform Relaxation Method for
  Time-Domain Analysis of Large Scale Integrated Circuits}, {IEEE Transactions
  on Computer-Aided Design of Integrated Circuits and Systems} 1 (1982)
  131--145.
\newblock \href {https://doi.org/10.1109/TCAD.1982.1270004}
  {\path{doi:10.1109/TCAD.1982.1270004}}.

\bibitem{Lumsdaine_White_WR_Semiconductor_NFMAO_1995}
A.~Lumsdaine, J.~White, {Accelerating Wave-Form relaxation methods with
  application to parallel semiconductor-device simulation}, {Numer. Funct.
  Anal. Optim.} {16}~({3-4}) ({1995}) {395--414}.
\newblock \href {https://doi.org/10.1080/01630569508816625}
  {\path{doi:10.1080/01630569508816625}}.

\bibitem{Miekkala_Nevanlinna_DI_CVGODE_1987}
U.~Miekkala, O.~Nevanlinna, {Convergence of Dynamic Iteration methods for
  initial-value problems}, {SIAM Journal on Scientific and Statistical
  Computing} {8}~({4}) ({1987}) {459--482}.
\newblock \href {https://doi.org/10.1137/0908046} {\path{doi:10.1137/0908046}}.

\bibitem{Miekkala_DI_LinearDAE_1989}
U.~Miekkala, {Dynamic Iteration methods applied to linear DAE systems}, {J.
  Comput. Appl. Math.} {25}~({2}) ({1989}) {133--151}.
\newblock \href {https://doi.org/10.1016/0377-0427(89)90044-7}
  {\path{doi:10.1016/0377-0427(89)90044-7}}.

\bibitem{Arnold_WR_2001}
M.~Arnold, {Constraint partitioning in dynamic iteration methods}, {Z. Angew.
  Math. Mech.} {81}~({3}) ({2001}) {S735--S738}.
\newblock \href {https://doi.org/10.1002/zamm.200108115143}
  {\path{doi:10.1002/zamm.200108115143}}.

\bibitem{Bartel_Brunk_Guenther_Shoeps_DI_ElectricCircuits_SIAMJSC_2013}
A.~Bartel, M.~Brunk, M.~Guenther, S.~Schoeps, {Dynamic Iteration for coupled
  problems of electrical circuits and Distributed Devices}, {SIAM J. Sci.
  Comput.} {35}~({2}) ({2013}) {B315--B335}.
\newblock \href {https://doi.org/10.1137/120867111}
  {\path{doi:10.1137/120867111}}.

\bibitem{Bartel_Guenther_PDAE_Electrical_SIAMreview_2018}
A.~Bartel, M.~Guenther, {PDAEs in Refined Electrical Network Modeling}, {SIAM
  Review} {60}~({1}) ({2018}) {56--91}.
\newblock \href {https://doi.org/10.1137/17M1113643}
  {\path{doi:10.1137/17M1113643}}.

\bibitem{Guenther_Bartel_Jacob_Reis_DI_PPHDAE_CTA_2021}
M.~Guenther, A.~Bartel, B.~Jacob, T.~Reis, {Dynamic iteration schemes and
  port-Hamiltonian formulation in coupled differential-algebraic equation
  circuit simulation}, {Int. J. Circuit Theory Appl.} {49}~({2}) ({2021})
  {430--452}.
\newblock \href {https://doi.org/10.1002/cta.2870}
  {\path{doi:10.1002/cta.2870}}.

\bibitem{Reichelt_White_Allen_WRsor_Parallel_SemiCond_SIAMJSC_1995}
M.~Reichelt, J.~White, J.~Allen, {Optimal convolution SOR acceleration of
  Wave-Form relaxation with application to parallelsimulation of
  semiconductor-devices }, {SIAM J. Sci. Comput.} {16}~({5}) ({1995})
  {1137--1158}.
\newblock \href {https://doi.org/10.1137/0916066} {\path{doi:10.1137/0916066}}.

\bibitem{Shoeps_DeGersem_Bartel_Cosimulation_I3ETMagnetics_2012}
S.~Schoeps, H.~De~Gersem, A.~Bartel, {Higher-Order Cosimulation of
  Field/Circuit Coupled Problems}, {IEEE Transactions in Magnetics} {48}~({2})
  ({2012}) {535--538}.
\newblock \href {https://doi.org/10.1109/TMAG.2011.2174039}
  {\path{doi:10.1109/TMAG.2011.2174039}}.

\bibitem{Jiang_Wing_WR_Spectra_LinearDAE_SIAMJNA_2000}
Y.~Jiang, O.~Wing, {A note on the spectra and pseudospectra of waveform
  relaxation operators for linear differential-algebraic equations}, {SIAM
  Journal on Numerical Analysis} {38}~({1}) ({2000}) {186--201}.
\newblock \href {https://doi.org/10.1137/S0036142997327063}
  {\path{doi:10.1137/S0036142997327063}}.

\bibitem{Jiang_WR_NonlinearDAE_I3ETCircuits_2004}
Y.~Jiang, {A general approach to Waveform Relaxation solutions of nonlinear
  Differential-Algebraic Equations: The continuous-time and discrete-time
  cases}, {IEEE Transaction on Circuits and Systems I} {51}~({9}) ({2004})
  {1770--1780}.
\newblock \href {https://doi.org/10.1109/TCSI.2004.834503}
  {\path{doi:10.1109/TCSI.2004.834503}}.

\bibitem{Lumsdaine_Wu_WR_Spectra_SIAMJSC_1997}
A.~Lumsdaine, D.~Wu, {Spectra and pseudospectra of waveform relaxation
  operators}, {SIAM J. Sci. Comput.} {18}~({1}) ({1997}) {286--304}.
\newblock \href {https://doi.org/10.1137/S106482759528778X}
  {\path{doi:10.1137/S106482759528778X}}.

\bibitem{Arnold_Gunther_preconDI_DAE_Bit_2001}
M.~Arnold, M.~Gunther, {Preconditioned dynamic iteration for coupled
  differential-algebraic systems}, {BIT} {41}~({1}) ({2001}) {1--25}.
\newblock \href {https://doi.org/10.1023/A:1021909032551}
  {\path{doi:10.1023/A:1021909032551}}.

\bibitem{Hout_WR_Nonlinear_AppNumMath_1995}
K.~Hout, {On the convergence of wave-form relaxation methods for stiff
  nonlinear ordinary differential-equations}, {Applied Numerical Mathematics}
  {18}~({1-3}) ({1995}) {175--190}.
\newblock \href {https://doi.org/10.1016/0168-9274(95)00052-V}
  {\path{doi:10.1016/0168-9274(95)00052-V}}.

\bibitem{Janssen_Vandewalle_WR_SOR_SIAMJNA_1997}
J.~Janssen, S.~Vandewalle, {On SOR waveform relaxation methods}, {SIAM Journal
  on Numerical Analysis} {34}~({6}) ({1997}) {2456--2481}.
\newblock \href {https://doi.org/10.1137/S0036142995294292}
  {\path{doi:10.1137/S0036142995294292}}.

\bibitem{Leimkuhler_WR_TimeStepAcceleration_SIAMJNA_1998}
B.~Leimkuhler, {Timestep acceleration of waveform relaxation}, {SIAM Journal on
  Numerical Analysis} {35}~({1}) ({1998}) {31--50}.
\newblock \href {https://doi.org/10.1137/S003614299528002X}
  {\path{doi:10.1137/S003614299528002X}}.

\bibitem{Lumsdaine_Wu_WR_KrylovAcceleration_SIAMJNA_2003}
A.~Lumsdaine, D.~Wu, {Krylov subspace acceleration of waveform relaxation},
  {SIAM Journal on Numerical Analysis} {41}~({1}) ({2003}) {90--111}.
\newblock \href {https://doi.org/10.1137/S0036142996313142}
  {\path{doi:10.1137/S0036142996313142}}.

\bibitem{Botchev_Oseledets_Tyrtyshnikov_WR_Krylov_ComputMathAppl_2014}
M.~A. Botchev, I.~V. Oseledets, E.~E. Tyrtyshnikov, {Iterative across-time
  solution of linear differential equations: Krylov subspace versus waveform
  relaxation}, {Comput. Math. Appl.} {67}~({12}) ({2014}) {2088--2098}.
\newblock \href {https://doi.org/10.1016/j.camwa.2014.03.002}
  {\path{doi:10.1016/j.camwa.2014.03.002}}.

\bibitem{Ladics_WR_SemilinearPDE_ErrorEstimate_ComputApplMath_2015}
T.~Ladics, {Error analysis of waveform relaxation method for semi-linear
  partial differential equations}, {J. Comput. Appl. Math.} {285} ({2015})
  {15--31}.
\newblock \href {https://doi.org/10.1016/j.cam.2015.02.003}
  {\path{doi:10.1016/j.cam.2015.02.003}}.

\bibitem{Bartel_Brunk_Schoeps_DI_cvgrate_JCAppMath_2014}
A.~Bartel, M.~Brunk, S.~Schoeps, {On the convergence rate of dynamic iteration
  for coupled problems with multiple subsystems}, {J. Comput. Appl. Math.}
  {262} ({2014}) {14--24}.
\newblock \href {https://doi.org/10.1016/j.cam.2013.07.031}
  {\path{doi:10.1016/j.cam.2013.07.031}}.

\bibitem{Ali_Bartel_Brunk_Schoeps_DI_cvg_ProcMI_2012}
G.~Ali, A.~Bartel, M.~Brunk, S.~Schoeps, {A Convergent Iteration Scheme for
  Semiconductor/Circuit Coupled Problems}, in: {Michielsen, B and Poirier, JR}
  (Ed.), {Scientific Computing in Electrical Engineering (SCEE 2010)},
  Vol.~{16} of {Mathematics in Industry-Cham}, {2012}, pp. {233--242}.
\newblock \href {https://doi.org/10.1007/978-3-642-22453-9\_\_25}
  {\path{doi:10.1007/978-3-642-22453-9\_\_25}}.

\bibitem{Gausling_Bartel_DI_COSIM_CG_Stoch_ProcMI_2018}
K.~Gausling, A.~Bartel, {Density Estimation Techniques in Cosimulation Using
  Spectral- and Kernel Methods}, in: {Langer, U and Amrhein, W and Zulehner, W}
  (Ed.), {Scientific Computing in Engineering, SCEE 2016}, Vol.~{28} of
  {Mathematics in Industry-Cham}, {2018}, pp. {81--89}.
\newblock \href {https://doi.org/10.1007/978-3-319-75538-0\_8}
  {\path{doi:10.1007/978-3-319-75538-0\_8}}.

\bibitem{Gausling_Bartel_DI_Cosim_cvgInterface_2016}
K.~Gausling, A.~Bartel, {Coupling Interfaces and Their Impact in Field/Circuit
  Co-Simulation}, {IEEE Trans. Magn.} {52}~({3}) ({MAR} {2016}).
\newblock \href {https://doi.org/10.1109/TMAG.2015.2471181}
  {\path{doi:10.1109/TMAG.2015.2471181}}.

\bibitem{Pade_Tischendorf_WR_CVGCriterion_NUMAlgo_2019}
J.~Pade, C.~Tischendorf, {Waveform relaxation: a convergence criterion for
  differential-algebraic equations}, {Numer. AlgorithmsNumer. Algorithms}
  {81}~({4, SI}) ({2019}) {1327--1342}.
\newblock \href {https://doi.org/10.1007/s11075-018-0645-5}
  {\path{doi:10.1007/s11075-018-0645-5}}.

\bibitem{tromeur_contrib_DTD}
M.~Garbey, D.~Tromeur-Dervout, On some {A}itken-like acceleration of the
  {S}chwarz method, Internat. J. Numer. Methods Fluids 40~(12) (2002)
  1493--1513.
\newblock \href {https://doi.org/10.1002/fld.407} {\path{doi:10.1002/fld.407}}.

\bibitem{tromeur_contrib_ISI:000268885300011}
D.~Tromeur-Dervout, {Meshfree Adaptative Aitken-Schwarz Domain Decomposition
  with application to Darcy Flow}, in: {Topping, BHV and Ivanyi, P} (Ed.),
  {Parallel, Distributed and Grid Computing for Engineering}, Vol.~{21} of
  {CSET Series}, Saxe-Coburg Publications, {2009}, pp. {217--250}.
\newblock \href {https://doi.org/10.4203/csets.21.11}
  {\path{doi:10.4203/csets.21.11}}.

\bibitem{tromeur_contrib_DTD2}
D.~Tromeur-Dervout, Approximating the trace of iterative solutions at the
  interfaces with nonuniform {F}ourier transform and singular value
  decomposition for cost-effectively accelerating the convergence of {S}chwarz
  domain decomposition, ESAIM: Proc. 42 (2013) 34--60.
\newblock \href {https://doi.org/10.1051/proc/201342004}
  {\path{doi:10.1051/proc/201342004}}.

\bibitem{Shourick_DD26}
H.~Shourick, D.~Tromeur-Dervout, L.~Chédot, {Aitken-Schwarz Heterogeneous
  Domain Decomposition for EMT-TS Simulation}, in: S.~Brenner, E.~Chung,
  A.~Klawonn, F.~Kwok, J.~Xu, J.~Zou (Eds.), Domain Decomposition Methods in
  Science and Engineering XXVI, Lecture Notes in Computational Sciences and
  Engineering, springer, 2022, to appear.

\bibitem{tromeur_contrib_RAS}
X.-C. Cai, M.~Sarkis, A restricted additive {S}chwarz preconditioner for
  general sparse linear systems, SIAM J. Sci. Comput. 21~(2) (1999) 792--797.
\newblock \href {https://doi.org/10.1137/S106482759732678X}
  {\path{doi:10.1137/S106482759732678X}}.

\bibitem{tromeur_contrib_Gander08schwarzmethods}
M.~J. Gander, \href{http://eudml.org/doc/130616}{Schwarz methods over the
  course of time}, Electronic Transactions on Numerical Analysis (2008)
  228--255.
\newline\urlprefix\url{http://eudml.org/doc/130616}

\end{thebibliography}

\end{document}